\documentclass[12pt]{amsart}

\usepackage{amsmath,amssymb,amsthm}
\usepackage{mathtools}
\usepackage{booktabs}
\usepackage{array}
\usepackage[hidelinks]{hyperref}
\usepackage{microtype}

\numberwithin{equation}{section}
\allowdisplaybreaks
\newtheorem{theorem}{Theorem}[section]
\newtheorem{proposition}[theorem]{Proposition}

\newtheorem{corollary}[theorem]{Corollary}
\theoremstyle{definition}
\newtheorem{definition}[theorem]{Definition}
\newtheorem{example}[theorem]{Example}
\newtheorem{remark}[theorem]{Remark}

\DeclareMathOperator{\codim}{codim}
\DeclareMathOperator{\tr}{tr}
\DeclareMathOperator{\Vol}{vol}

\newcommand{\R}{\mathbb{R}}
\newcommand{\C}{\mathbb{C}}
\newcommand{\Z}{\mathbb{Z}}
\newcommand{\HH}{\mathbb{H}}
\newcommand{\cM}{\mathcal{M}}
\newcommand{\cN}{\mathcal{N}}
\newcommand{\cS}{\mathcal{S}}
\newcommand{\fkg}{\mathfrak{g}}
\newcommand{\fkk}{\mathfrak{k}}
\newcommand{\fkp}{\mathfrak{p}}
\newcommand{\Skew}{\mathrm{Skew}}
\newcommand{\St}{\mathrm{St}}
\newcommand{\Gr}{\mathrm{Gr}}
\newcommand{\id}{\mathrm{id}}
\makeatletter
\def\author@andify{%
  \nxandlist{\unskip ,\penalty-1 \space\ignorespaces}%
    {\unskip ,\penalty-2 \space\ignorespaces}%
    {\unskip ,\penalty-2 \space\ignorespaces}%
}
\makeatother

\begin{document}

\title{Almost-Global Coordinates on Matrix Manifolds}

\author{Bingjie Li}
\address{School of Public Health, Tongji Medical College\newline
Huazhong University of Science and Technology}
\email{bingjieli@hust.edu.cn}

\author{Jiaxiang Wang}
\address{School of Mathematical Sciences and LPMC, Nankai University\newline
Tianjin 300071, People's Republic of China}
\email{wangjx\_manifold@126.com}

\author{Zhigang Yao}
\address{Department of Statistics and Data Science\newline
National University of Singapore}
\email{zhigang.yao@nus.edu.sg}


\begin{abstract}
We study global coordinates on smooth manifolds outside some closed subsets of measure zero. We first give an elementary proof that every compact smooth manifold admits such coordinates, allowing the Euclidean parameter domain to be disconnected. We then study parametrizations of classical matrix groups, Stiefel manifolds, Grassmannians and real flag manifolds. For these parametrizations, we describe the inverse maps and compute the codimensions of the omitted subsets. We construct these parametrizations using Cayley transforms,
graphs of linear maps and block Gaussian elimination.
We also give formulas for Riemannian gradients and Haar integrals
in these coordinates.
\end{abstract}

\maketitle
\tableofcontents
\clearpage

\section{Introduction}
\label{sec:introduction}

A nonempty compact manifold generally has no global
coordinate chart. As a model case, stereographic projection identifies
$S^d$ minus one point with $\R^d$. This suggests studying coordinates that cover a manifold except for a closed subset of measure zero. These coordinates can be used for integration and for continuous optimization on a manifold. However, critical points may lie in the omitted subset. We therefore ask whether such coordinates exist, whether they can be written explicitly, and how their inverses and omitted subsets can be described. 

Let $\cM$ be a smooth $d$-manifold. We call a diffeomorphism
$$
\mathbf G:\cS\longrightarrow\cM\setminus\Sigma
$$
a \emph{strict approximate parametrization} if $\Sigma$ is a closed
subset of $\cM$ of Lebesgue measure zero and $\cS$ is an open subset
of $\R^d$.
The parameter domain $\cS$ may be disconnected, and even when connected,
it need not be contractible. The term ``approximate'' refers only to
the omitted set; both the map and its inverse are exact. We also use
a broader notion that allows noninjective local diffeomorphisms,
but the explicit parametrizations below are strict. 

\subsection{Existence and the choice of the exceptional set}

Before presenting the explicit constructions, we state an existence theorem for strict approximate parametrizations
of smooth manifolds.

\begin{theorem}
\label{thm:existence}
Let $\cM$ be a compact smooth manifold of dimension $d\geq1$.
There exist an open set $\cS\subset\R^d$ and a diffeomorphism
$\mathbf G:\cS\to\cM\setminus\Sigma$, where $\Sigma$ is a finite
union of compact embedded $(d-1)$-dimensional submanifolds of $\cM$.
In particular, $\Sigma$ is closed and $\dim_H\Sigma\leq d-1$.
\end{theorem} 

The existence of strict approximate parametrizations follows
from the classical cut-locus construction. Every smooth manifold
admits a complete Riemannian metric. On a complete connected
Riemannian manifold, the exponential map at a point $p$ restricts
to a diffeomorphism from its open injectivity domain onto the
complement of $\operatorname{Cut}(p)$. The cut locus is closed
and satisfies $\dim_H\operatorname{Cut}(p)\leq d-1$;
see~\cite[Theorem~10.34(a),(c), p.~311]{Lee18} and~\cite[Lemma~3]{IT98}.
For a disconnected manifold, we apply this construction to each
connected component.

The codimension of the exceptional set depends on the chosen
parametrization. On $S^d$, stereographic projection omits a point,
whereas parametrizing the two open hemispheres omits an equator.
The corresponding codimensions are $d$ and $1$, respectively.
For the matrix examples, we determine
codimension from the dimensions of the smooth part of the
omitted set. 

\subsection{Explicit coordinates on matrix groups}

The first construction is the Cayley transform. We write $A^*$ for
the transpose over $\R$ and the conjugate transpose over $\C$ or $\HH$.
For an $n\times n$ skew-adjoint matrix $X$ (that is, $X^*=-X$), define
$$
\mathbf C_n(X)=(I+X)^{-1}(I-X),
$$
where $I$ denotes the identity matrix.

If $Q=\mathbf C_n(X)$, then
$$
X=\mathbf C_n^{-1}(Q)=(I+Q)^{-1}(I-Q).
$$
The following theorem shows that $\mathbf C_n$ gives strict
approximate parametrizations of $SO(n)$, $U(n)$ and $Sp(n)$.
In each case, the omitted set consists precisely of matrices $Q$
for which $I+Q$ is singular, or equivalently, for which $-1$
is an eigenvalue.

\begin{theorem}
\label{thm:main-classical}
Let $G$ be one of the groups $SO(n)$, $U(n)$ and $Sp(n)$,
and let $\fkg$ denote its Lie algebra of skew-adjoint matrices.
Set
$$
G_+=\{Q\in G:I+Q\text{ is invertible}\}.
$$
The Cayley transform is a real-analytic diffeomorphism
$\mathbf C_n:\fkg\to G_+$, with inverse
$$
\mathbf C_n^{-1}(Q)=(I+Q)^{-1}(I-Q).
$$
The complement $G\setminus G_+$ is closed and has measure zero.
Its real codimension is $1$ for $SO(n)$ with $n\geq2$ and
for $U(n)$ with $n\geq1$, and is $3$ for $Sp(n)$ with $n\geq1$.
For $n\geq2$, the Cayley chart on $SO(n)$ and its translate
by a reflection together give a strict approximate
parametrization of $O(n)$.
\end{theorem}

To compute the codimension, we prescribe the dimension of the
$-1$ eigenspace, choose this eigenspace, and then choose the
orthogonal, unitary or quaternionic unitary transformation on
its orthogonal complement. Counting these parameters gives
the dimension of each stratum. For $SO(n)$, the $-1$ eigenspace
has even dimension. For $O(n)$, a separate translate is needed
to cover the component of determinant $-1$. We prove these
assertions in Section~\ref{sec:lie-groups}. Combining the
orthogonal parametrization with the QR decomposition with
positive diagonal entries also gives coordinates on $GL(n,\R)$
and $SL(n,\R)$.

We next consider the special unitary group $SU(n)$. If
$X\in\mathfrak{u}(n)$ has eigenvalues $i\lambda_1,\ldots,i\lambda_n$,
where $\lambda_j\in\R$, then
$$
\det\mathbf C_n(X)
=\exp\left(-2i\sum_{j=1}^n\arctan\lambda_j\right).
$$
Here and below, $\arctan$ denotes the principal branch with values
in $(-\pi/2,\pi/2)$. For Hermitian matrices, it is defined by
applying this function to the eigenvalues in a unitary
diagonalization. For $n\geq3$, the condition $\tr X=0$
does not imply $\det\mathbf C_n(X)=1$. To impose the
determinant-one condition, we apply the Cayley transform to
$Y+isI$, where $Y\in\mathfrak{su}(n)$ and $s\in\R$ is chosen
so that $\det\mathbf C_n(Y+isI)=1$.

\begin{theorem}
\label{thm:main-special-unitary}
Let $n\geq2$, and set
$$
K_n=\{k\in\Z:|k|<n/2\},\qquad
SU_+(n)=\{Q\in SU(n):I+Q\text{ is invertible}\}.
$$
For every $Y\in\mathfrak{su}(n)$ and $k\in K_n$, there exists
a unique $s_k(Y)\in\R$ such that
$$
\sum_{j=1}^n\arctan(\mu_j+s_k(Y))=k\pi,
$$
where $\mu_1,\ldots,\mu_n$ are the eigenvalues of the Hermitian
matrix $-iY$, counted with multiplicity, and $\arctan$ takes
values in $(-\pi/2,\pi/2)$.

Each function $s_k:\mathfrak{su}(n)\to\R$ is real-analytic.
The map
$$
\mathbf C_n^{SU}:K_n\times\mathfrak{su}(n)\longrightarrow SU_+(n),
\qquad
\mathbf C_n^{SU}(k,Y)=\mathbf C_n(Y+i s_k(Y)I),
$$
is a real-analytic diffeomorphism. For each fixed $k\in K_n$,
the map $Y\mapsto\mathbf C_n^{SU}(k,Y)$ maps
$\mathfrak{su}(n)$ diffeomorphically onto one connected
component of $SU_+(n)$. As $k$ ranges over $K_n$, these images
give all the connected components, without repetition.

The omitted set $SU(n)\setminus SU_+(n)$ has real codimension
$3$ for $n=2$ and $1$ for $n\geq3$.
\end{theorem}

We first explain the existence and uniqueness of $s_k(Y)$.
For fixed $Y$, let $\mu_1,\ldots,\mu_n$ be the eigenvalues
of $-iY$. Then
$$
\frac{d}{ds}\sum_{j=1}^n\arctan(\mu_j+s)
=\sum_{j=1}^n\frac{1}{1+(\mu_j+s)^2}>0,
$$
and
$$
\lim_{s\to\pm\infty}\sum_{j=1}^n\arctan(\mu_j+s)
=\pm\frac{n\pi}{2}.
$$
Hence the sum is strictly increasing in $s$ and has range
$(-n\pi/2,n\pi/2)$. It therefore equals $k\pi$ at a unique
value $s=s_k(Y)$ precisely when $k\in K_n$. The determinant
identity above then gives
$\det\mathbf C_n(Y+i s_k(Y)I)=1$.

The sum can also be written as $\tr\arctan(-iY+sI)$,
which is real-analytic in the matrix entries of $Y$ and in $s$.
Since its derivative with respect to $s$ is positive, the
real-analytic implicit function theorem gives the real-analytic
dependence of $s_k(Y)$ on $Y$.

We next describe the inverse map. Given $Q\in SU_+(n)$, set
$$
\begin{array}{l}
X=(I+Q)^{-1}(I-Q),\\[4pt]
s=\displaystyle\frac{\tr X}{in},\qquad Y=X-isI,\\[6pt]
k=\displaystyle\frac1\pi\tr\arctan(-iX).
\end{array}
$$
The ordinary Cayley inverse gives $X\in\mathfrak{u}(n)$,
so $s\in\R$ and $Y\in\mathfrak{su}(n)$. Since $\det Q=1$,
the determinant identity implies that $k$ is an integer.
The range of $\arctan$ gives $|k|<n/2$, and hence $k\in K_n$.
By uniqueness in the scalar equation, $s=s_k(Y)$.
Thus these formulas recover the unique pair $(k,Y)$
corresponding to $Q$, without choosing eigenvectors.

The integer $k$ is locally constant on $SU_+(n)$, and for
each fixed $k$ the parameter space $\mathfrak{su}(n)$ is
connected. Together with the inverse formulas, this shows
that the images for distinct $k$ are precisely the connected
components of $SU_+(n)$. The full proof, including the
codimension calculation, is given in
Theorem~\ref{thm:SUn-cayley}.

\subsection{Stiefel manifolds, Grassmannians and flag
manifolds}

We next consider Stiefel manifolds, Grassmannians and flag
manifolds. These spaces can be described as quotients of matrix
groups. However, a set of measure zero in a group may project
onto the entire quotient, so the group coordinates above do
not automatically give strict approximate parametrizations
of these spaces. We therefore construct each parametrization
and its inverse directly.

Let $\mathbb F=\R,\C$ or $\HH$. We denote by
$\St_{\mathbb F}(n,p)$ the manifold of orthonormal $p$-frames
in $\mathbb F^n$, and by $\Gr_{\mathbb F}(n,p)$ the manifold
of $p$-dimensional $\mathbb F$-linear subspaces of $\mathbb F^n$.
Set $\beta=\dim_\R\mathbb F$, so that $\beta=1,2$ or $4$,
respectively. 

\begin{theorem}
\label{thm:main-homogeneous}
Let $1\leq p<n$.
\begin{enumerate}
\item[(i)] For $A\in\mathbb F^{p\times p}$ with $A^*=-A$
and $B\in\mathbb F^{(n-p)\times p}$, the map
$$
(A,B)\longmapsto
\begin{pmatrix}
2(I_p+A+B^*B)^{-1}-I_p\\
2B(I_p+A+B^*B)^{-1}
\end{pmatrix}
$$
is a strict approximate parametrization of
$\St_{\mathbb F}(n,p)$. Its image consists of the frames
$$
Q=\begin{pmatrix}Q_1\\Q_2\end{pmatrix},
\qquad Q_1\in\mathbb F^{p\times p},
$$
for which $I_p+Q_1$ is invertible.

\item[(ii)] The graph map
$$
B\longmapsto\{(v,Bv):v\in\mathbb F^p\},
\qquad B\in\mathbb F^{(n-p)\times p},
$$
is a strict approximate parametrization of
$\Gr_{\mathbb F}(n,p)$. Its image consists of the subspaces
transverse to $\{0\}\times\mathbb F^{n-p}$.
\end{enumerate}
The omitted subsets have the following real codimensions:
$$
\begin{array}{c@{\quad}|@{\quad}c@{\quad}c@{\quad}c}
&\mathbb F=\R&\mathbb F=\C&\mathbb F=\HH\\ \hline
\St_{\mathbb F}(n,p)&n-p&2(n-p)+1&4(n-p)+3\\
\Gr_{\mathbb F}(n,p)&1&2&4
\end{array}
$$
\end{theorem}

The inverse maps and the codimension calculations are given
in Section~\ref{sec:symmetric-spaces}. We also construct
strict approximate parametrizations of nontrivial real
partial flag manifolds and the real Lagrangian Grassmannian
$U(n)/O(n)$. In both cases, the omitted subsets have real
codimension one.

\subsection{Relation to earlier work and applications}

The explicit constructions have several classical antecedents.
Cayley~\cite{Cay46} introduced the rational orthogonal formula.
Miller~\cite[Theorem~A]{Miller85} describes the successive strata
of his filtration of Stiefel manifolds as vector bundles over
Grassmannians. Schwede~\cite[Propositions~B.8 and B.11]{Sch22}
gives explicit Cayley maps and describes the kernel filtration
uniformly over $\R$, $\C$ and $\HH$. The Stiefel formulas used
here agree with Schwede's maps up to a fixed sign change in the
target and a rescaling of the parameters. Mac\'ias-Virg\'os,
Pereira-S\'aez and Tanr\'e~\cite{MVPT18} study Cayley charts
on Stiefel manifolds over these three fields, with a specified
domain of injectivity at each base point. Kume and
Yamada~\cite{KY24} develop localized Cayley charts and explicit
inverse maps for real Stiefel manifolds.

The Grassmann and real flag graph charts are open Schubert cells;
see Fomenko and Fuchs~\cite[Sections~5.4C--D]{FF16}.
The symmetric graph chart on $U(n)/O(n)$ and its intersection
strata are described by Arnol'd~\cite[Corollary~3.1.3 and
Lemma~3.2.1]{Arn67}. For $SU(3)$, Sch\"afers, Peardon and
G\"unther~\cite{SPG25} use a phase modification of the Cayley
transform for molecular dynamics simulations. Their transform gives
a local parametrization near the identity. Our scalar-shift
construction describes all branches on the open set of matrices
in $SU(n)$ without eigenvalue $-1$.

These coordinates also have direct computational applications.
Wen and Yin~\cite{WY13} use Cayley updates to preserve
orthogonality in optimization. Jauch, Hoff and Dunson~\cite{JHD20}
use injective Cayley parametrizations with null complement and
their Jacobians to simulate distributions on real Stiefel and
Grassmann manifolds. T\^oyama~\cite[Theorems~1--3]{Toy48}
computes Haar densities, up to normalization, in Cayley parameters
for unitary, compact symplectic and orthogonal groups.
Manifold fitting provides a statistical setting in which local
geometric estimation can introduce Grassmann and Stiefel
variables~\cite{YSLY23,YX25}. Section~\ref{subsec:opt}
illustrates how the Grassmann coordinates developed here represent
a smooth local subspace-fitting objective in Euclidean parameters.

Our purpose is to give a common, explicit treatment of these
coordinates, their inverses and their omitted sets. We retain
the component restrictions and endpoint cases, and distinguish
measure-theoretic conclusions from differential estimates. In
particular, growth of the inverse differential does not by itself
imply growth of the condition number. In the applications, we
compute the pullback metric and the Cayley singular values before
using them for gradients and Haar integration.

Section~\ref{sec:theory} contains the definitions, the proof of
Theorem~\ref{thm:existence} and the integration and optimization
facts used later. Section~\ref{sec:compatibility} treats products,
quotients and bundles. Section~\ref{sec:cayley} establishes the
common Cayley identities. Sections~\ref{sec:lie-groups}
and~\ref{sec:symmetric-spaces} give the explicit group and
homogeneous-space constructions. The applications are discussed in Section~\ref{sec:discussion}.

\section{Coordinates outside some measure zero subsets}
\label{sec:prelim}
\label{sec:theory}

In this section, we prove Theorem~\ref{thm:existence} and record
the facts on integration and optimization used later.

\subsection{Measure and dimension}
\label{subsec:manifolds}

All manifolds are Hausdorff, second countable, smooth and without
boundary. Unless stated otherwise, their dimension is positive.

Measure zero on a manifold is understood in local coordinates.
Every Riemannian volume density is locally a smooth positive
multiple of Lebesgue measure, so this notion is independent of
the choice of metric and requires no orientability assumption.
Diffeomorphisms preserve sets of measure zero, since they and
their inverses are locally Lipschitz in coordinates.

\begin{proposition}
\label{prop:submanifold-measure}
Every embedded submanifold of dimension less than $d$ has measure
zero in a $d$-manifold. The same holds for a countable union of
such submanifolds.
\end{proposition}

\begin{proof}
In local coordinates, a $k$-dimensional embedded submanifold
lies in $\R^k\times\{0\}$, which has Lebesgue measure zero
for $k<d$. The assertion follows from second countability
and countable subadditivity.
\end{proof}

When an omitted set is a finite union of smooth strata, its
dimension is the largest stratum dimension, and its codimension
is the ambient dimension minus this number. For more general
omitted sets, we state bounds for their Hausdorff dimension.
We assign codimension $+\infty$ to the empty set.

\subsection{Approximate parametrizations}
\label{subsec:definitions}

\begin{definition}
\label{def:approx-param}
Let $\cM$ be a smooth $d$-manifold, $\cS\subset\R^d$ an open
set, and $\Sigma\subset\cM$ a closed set of measure zero.
A diffeomorphism
$$
\mathbf G:\cS\longrightarrow\cM\setminus\Sigma
$$
is called a \emph{strict approximate parametrization} of $\cM$.
More generally, a smooth map $\mathbf G:\cS\to\cM$ is called
an \emph{approximate parametrization} with exceptional set
$\Sigma$ if the restriction
$$
\mathbf G:\mathbf G^{-1}(\cM\setminus\Sigma)
\longrightarrow\cM\setminus\Sigma
$$
is a surjective local diffeomorphism. It is \emph{finite-to-one}
if every point of $\cM\setminus\Sigma$ has finitely many preimages.
\end{definition}

The parameter domain $\cS$ may be disconnected, and even when
connected, it need not be contractible. We also allow parameter
spaces diffeomorphic to open subsets of $\R^d$.

Since $\Sigma$ is closed and has measure zero,
$\cM\setminus\Sigma$ is open and dense in $\cM$.
For a strict approximate parametrization, $\Sigma$ is exactly
the omitted set, and both the map and its inverse are smooth.
Thus ``approximate'' refers only to the omitted set.
All explicit parametrizations below are strict.

\subsection{Existence on compact manifolds}
\label{subsec:existence}

Theorem~\ref{thm:existence} follows by making a finite family
of coordinate pieces disjoint. We remove only points on the
boundaries of coordinate balls.

\begin{proof}[Proof of Theorem~\ref{thm:existence}]
Choose coordinate balls $U_1,\ldots,U_N$ covering $\cM$.
For each $i$, choose a coordinate map $\varphi_i$ defined on
a neighborhood of $\overline{U_i}$ with
$\varphi_i(U_i)=B(0,1)$ and
$\varphi_i(\partial U_i)=\partial B(0,1)$. Set
$$
\Sigma=\bigcup_{i=1}^N\partial U_i,\qquad
V_i=U_i\setminus\left(\Sigma\cup\bigcup_{j<i}\overline{U_j}\right).
$$
Each $V_i$ is open, and the $V_i$ are pairwise disjoint.
If $x\notin\Sigma$, take the least $i$ such that $x\in U_i$.
Then $x\notin\overline{U_j}$ for $j<i$, because it belongs
to neither $U_j$ nor $\partial U_j$. Hence
$$
\cM\setminus\Sigma=\coprod_{i=1}^N V_i.
$$

Choose vectors $a_i$ so that the balls $a_i+B(0,1)$ are disjoint.
On the open set
$$
\cS=\bigcup_{i=1}^N\bigl(a_i+\varphi_i(V_i)\bigr),
$$
define $\mathbf G(a_i+z)=\varphi_i^{-1}(z)$.
The inverse on $V_i$ is $x\mapsto a_i+\varphi_i(x)$.
Both maps are smooth on their disjoint open pieces, so they give
a diffeomorphism $\cS\cong\cM\setminus\Sigma$.

Each $\partial U_i$ is a compact embedded hypersurface. It is
null by Proposition~\ref{prop:submanifold-measure} and has
Hausdorff dimension $d-1$. Their finite union is closed and
null and has Hausdorff dimension at most $d-1$.
\end{proof}

\begin{remark}
\label{rem:existence-generic}
The construction gives an open subset of Euclidean space,
possibly disconnected. The boundaries need not meet transversely,
so we state a Hausdorff dimension bound without asserting that there is a smooth
stratification of their union. The condition $d\geq1$ is also
needed: a nonempty open subset of $\R^0$ is one point, while
a compact zero-dimensional manifold may have several points.
Its volume is counting measure, and no nonempty set is null.
\end{remark}

\subsection{Local inverses and critical values}
\label{subsec:sard}

We recall Sard's theorem on the critical values of smooth maps;
see~\cite[Section~2, pp.~10--11]{Mil65}.

\begin{theorem}[Sard's theorem]
\label{thm:sard}
Let $F:\cM\to\cN$ be a smooth map. Then the set of critical
values of $F$ has measure zero in $\cN$.
\end{theorem}

\begin{corollary}
\label{cor:sard}
The regular values of $F$ have complement of measure zero in
$\cN$. If $F$ is proper, they form an open set, and the
restriction of $F$ over this set is a submersion.
\end{corollary}

\begin{proof}
The first assertion follows from Sard's theorem. If $F$ is
proper, the image of its closed critical set is closed. Hence
the set of regular values is open. On its inverse image, the
differential of $F$ is surjective, so the restriction is a
submersion.
\end{proof}

\begin{remark}
\label{rem:sard-regularity}
When $\cM$ and $\cN$ have the same dimension, the inverse
function theorem gives a local inverse at every preimage of a
regular value. This does not imply injectivity. For example,
the map $\R\to S^1$ given by $t\mapsto e^{it}$ is a local
diffeomorphism, but every point has infinitely many preimages.
For the strict approximate parametrizations below, we construct
the global inverses explicitly and verify that the omitted
sets have measure zero.
\end{remark}

\subsection{Metrics, integration and optimization}

We next explain how approximate parametrizations can be used
for integration and optimization. We first describe the
Riemannian gradient in the parameter domain.

\begin{proposition}
\label{prop:pullback}
Let $\mathbf G:\cS\to\cM$ be an approximate parametrization
with exceptional set $\Sigma$, and let $g$ be a Riemannian
metric on $\cM$. Then $\mathbf G^*g$ is a Riemannian metric
on $\mathbf G^{-1}(\cM\setminus\Sigma)$. For every
$f\in C^\infty(\cM)$ and
$x\in\mathbf G^{-1}(\cM\setminus\Sigma)$, we have
$$
\nabla_{\mathbf G^*g}(f\circ\mathbf G)(x)
=(d\mathbf G_x)^{-1}\bigl(\nabla_g f(\mathbf G(x))\bigr).
$$
The inverse of $d\mathbf G_x$ is also its adjoint with respect
to $\mathbf G^*g$ and $g$.
\end{proposition}

\begin{proof}
For $x\in\mathbf G^{-1}(\cM\setminus\Sigma)$, the differential
$d\mathbf G_x$ is invertible. Hence $\mathbf G^*g$ is positive
definite there. For any $v\in T_x\cS$, the chain rule gives
$$
\begin{aligned}
d(f\circ\mathbf G)_x(v)
&=g_{\mathbf G(x)}\bigl(\nabla_g f(\mathbf G(x)),d\mathbf G_xv\bigr)\\
&=(\mathbf G^*g)_x\bigl((d\mathbf G_x)^{-1}
\nabla_g f(\mathbf G(x)),v\bigr).
\end{aligned}
$$
The gradient formula follows from the definition. Since
$d\mathbf G_x$ is an isometry for these two metrics, its
inverse equals its adjoint.
\end{proof}

\begin{theorem}
\label{thm:measure-equiv}
Let $\mathbf G:\cS\to\cM\setminus\Sigma$ be a strict
approximate parametrization. Equip $\cS$ and $\cM$ with
Riemannian metrics, and denote the corresponding volume
measures by $\mu_{\cS}$ and $\mu_{\cM}$. Let
$J_{\mathbf G}=|\det d\mathbf G|$, computed in orthonormal
frames for these metrics. Then the following hold.
\begin{enumerate}
\item[(i)] For every $\mu_{\cM}$-integrable function $f$,
or nonnegative measurable function $f$,
$$
\int_{\cM}f\,d\mu_{\cM}
=\int_{\cS}(f\circ\mathbf G)J_{\mathbf G}\,d\mu_{\cS}.
$$
For nonnegative $f$, both sides may be infinite.
\item[(ii)] For every continuous function $f:\cM\to\R$,
$$
\inf_{\cS}(f\circ\mathbf G)=\inf_{\cM}f.
$$
Consequently, a sequence $(x_j)\subset\cS$ is minimizing
for $f\circ\mathbf G$ if and only if $(\mathbf G(x_j))$
is minimizing for $f$.
\end{enumerate}
\end{theorem}

\begin{proof}
For (i), apply the change-of-variables formula to
$\mathbf G:\cS\to\cM\setminus\Sigma$. Since $\Sigma$ has
measure zero, its removal does not change the integral over
$\cM$. This proves the formula for nonnegative measurable
functions. Applying it to positive and negative parts, and
then to real and imaginary parts, gives the integrable case.
The use of volume densities requires no orientability assumption.

For (ii), note that $\cM\setminus\Sigma$ is dense in $\cM$.
For each $p\in\cM$, choose $p_j\in\cM\setminus\Sigma$
with $p_j\to p$. By continuity,
$$
\inf_{\cS}(f\circ\mathbf G)
=\inf_{\cM\setminus\Sigma}f
\leq\lim_{j\to\infty}f(p_j)=f(p).
$$
Taking the infimum over $p\in\cM$ and using the reverse
inequality gives (ii), also when $\inf_{\cM}f=-\infty$.
The assertion about minimizing sequences follows from this
equality.
\end{proof}

Without continuity, the infimum may change after removing
$\Sigma$. For example, a function equal to $-1$ at the point
omitted by stereographic projection and $0$ elsewhere has
infimum $-1$, while its pullback is identically zero.
For a finite-to-one parametrization, the integration formula
must account for the number of preimages of each point.
Even for a continuous function, equality of infima does not
imply that the infimum is attained in $\cS$: a minimizer on
$\cM$ may lie in $\Sigma$.

\section{General constructions of approximate parametrizations}
\label{sec:compatibility}
\label{subsec:algebra}

In this section, we construct approximate parametrizations of
products, quotients and bundles. For a quotient, we require the
image of the given parametrization to be invariant under the
group action. For a bundle, we require a trivialization over
the complement of the exceptional set in the base.

\subsection{Products}

\begin{theorem}
\label{thm:product}
Let $\mathbf G_i:\cS_i\to\cM_i$ be approximate
parametrizations with exceptional sets $\Sigma_i$, $i=1,2$.
Then
$$
\mathbf G_1\times\mathbf G_2:\cS_1\times\cS_2
\longrightarrow\cM_1\times\cM_2
$$
is an approximate parametrization with exceptional set
$$
\Sigma_{12}=(\Sigma_1\times\cM_2)
\cup(\cM_1\times\Sigma_2).
$$
It is strict if both factors are strict, and finite-to-one if
both factors are finite-to-one. If each nonempty $\Sigma_i$ is
a finite union of smooth strata, then
$$
\codim\Sigma_{12}=\min\{\codim\Sigma_1,\codim\Sigma_2\}.
$$
Here the empty set has codimension $+\infty$, as in
Section~\ref{sec:theory}.
\end{theorem}

\begin{proof}
The complement of $\Sigma_{12}$ is the open set
$(\cM_1\setminus\Sigma_1)\times(\cM_2\setminus\Sigma_2)$,
and Fubini's theorem in local coordinates shows that
$\Sigma_{12}$ has measure zero. The product map is onto this
complement, and its differential is the direct sum of the two
invertible differentials. If both maps are strict, the product
of their inverses gives a global inverse. If both maps are
finite-to-one, every fiber over the complement of $\Sigma_{12}$
is a product of two finite sets.

A stratum of dimension $r$ in $\Sigma_1$ gives a subset of
dimension $r+\dim\cM_2$ in $\Sigma_1\times\cM_2$.
Applying the same observation to $\cM_1\times\Sigma_2$
and taking the largest dimension proves the codimension formula.
\end{proof}

\subsection{Quotients}

Let $H$ be a closed subgroup of a Lie group $G$, and let
$\pi:G\to G/H$ be the quotient map. This is a locally trivial
principal $H$-bundle; see~\cite[Section~I.10]{Kna23}.
A local section $s:V\to G$ gives the trivialization
$(q,h)\mapsto s(q)h$ of $\pi^{-1}(V)$.

\begin{theorem}
\label{thm:quotient}
Let $\mathbf G:\cS\to U=G\setminus\Sigma$ be a strict
approximate parametrization of $G$. Assume that $U$ is
invariant under right multiplication by $H$. Then
$$
x\cdot h=\mathbf G^{-1}(\mathbf G(x)h)
$$
defines a smooth free proper $H$-action on $\cS$, and
$$
\overline{\mathbf G}:\cS/H\longrightarrow
(G/H)\setminus\pi(\Sigma),\qquad
[x]\longmapsto\pi(\mathbf G(x))
$$
is a diffeomorphism. The set $\pi(\Sigma)$ is closed and has
measure zero in $G/H$. If $\psi:T\to\cS/H$ is a
diffeomorphism from an open set
$T\subset\R^{\dim G-\dim H}$, then
$\overline{\mathbf G}\circ\psi$ is a strict approximate
parametrization of $G/H$ with exceptional set $\pi(\Sigma)$.
\end{theorem}

\begin{proof}
The action on $\cS$ is obtained from right multiplication on
$U$ by the diffeomorphism $\mathbf G$. It is smooth and free.
Since $H$ is closed, the map $(g,h)\mapsto(g,gh)$ identifies
$U\times H$ with the closed subset
$$
\{(g,g')\in U\times U:g^{-1}g'\in H\}.
$$
Hence the action is proper. Local bundle trivializations give
a diffeomorphism $U/H\to\pi(U)$, and therefore give the
asserted map $\overline{\mathbf G}$.

Both $U$ and $\Sigma$ are unions of fibers of $\pi$. Thus
$$
\pi(\Sigma)=(G/H)\setminus\pi(U).
$$
Since $\pi$ is open, $\pi(\Sigma)$ is closed. In a local
trivialization over $V\subset G/H$, its preimage is
$(\pi(\Sigma)\cap V)\times H$ and has measure zero.
Restricting the second factor to a coordinate ball of positive
finite volume, Fubini's theorem implies that
$\pi(\Sigma)\cap V$ has measure zero. If $H$ is discrete,
we use a single point with counting measure. A countable
family of trivializations proves that $\pi(\Sigma)$ has
measure zero in $G/H$. The last assertion now follows by
composition with $\psi$.
\end{proof}

\begin{remark}
\label{rem:quotient-explicit}
The space $\cS/H$ need not be diffeomorphic to a Euclidean
open set, so this condition is needed to obtain a strict
approximate parametrization. Invariance of $U$ is also needed
in the proof of the measure-zero assertion. For example, under
the quotient of $\R^2$ by the vertical subgroup
$\{0\}\times\R$, the horizontal axis has measure zero but
projects onto the whole quotient. For the homogeneous spaces
in Section~\ref{sec:symmetric-spaces}, we construct the
parametrizations and their inverses directly.
\end{remark}

\subsection{Bundles}

\begin{theorem}
\label{thm:fiberbundle}
Let $\pi:E\to B$ be a smooth bundle with fiber $F$, and let
$\mathbf G_B:\cS_B\to B$ and $\mathbf G_F:\cS_F\to F$
be approximate parametrizations with exceptional sets
$\Sigma_B$ and $\Sigma_F$. Assume that there is a smooth
trivialization
$$
\tau:(B\setminus\Sigma_B)\times F
\longrightarrow\pi^{-1}(B\setminus\Sigma_B).
$$
Set $\cS_B^\circ=\mathbf G_B^{-1}(B\setminus\Sigma_B)$.
Then
$$
\mathbf G_E:\cS_B^\circ\times\cS_F\longrightarrow E,
\qquad
\mathbf G_E(x,y)=\tau(\mathbf G_B(x),\mathbf G_F(y))
$$
is an approximate parametrization with exceptional set
$$
\Sigma_E=\pi^{-1}(\Sigma_B)
\cup\tau((B\setminus\Sigma_B)\times\Sigma_F).
$$
It is strict if both given parametrizations are strict.
\end{theorem}

\begin{proof}
The complement of $\Sigma_E$ is the open set
$\tau((B\setminus\Sigma_B)\times(F\setminus\Sigma_F))$,
so $\Sigma_E$ is closed. In local bundle coordinates,
$\pi^{-1}(\Sigma_B)$ is a product of a measure-zero set
with the fiber. Hence it has measure zero by Fubini's theorem
and a countable cover. The second part of $\Sigma_E$ has
measure zero because $\tau$ is a diffeomorphism.

Over $E\setminus\Sigma_E$, the map $\mathbf G_E$ is
surjective, and its differential is $d\tau$ composed with the
direct sum of the two invertible differentials. Thus it is an
approximate parametrization. If both given maps are strict,
then $\cS_B^\circ=\cS_B$, and
$(\mathbf G_B^{-1}\times\mathbf G_F^{-1})\circ\tau^{-1}$
is its global inverse.
\end{proof}

\section{The Cayley transform and its inverse}
\label{sec:cayley}
\label{subsec:cayley-general}

We first prove the matrix identities used in the group constructions.
The same calculation applies over $\R$, $\C$ and $\HH$.
In the next section, we describe the omitted sets and compute
their codimensions.

\subsection{Matrix conventions}
\label{subsec:lie}

Let $\mathbb F=\R$, $\C$ or $\HH$. As in the Introduction,
$A^*$ denotes the transpose over $\R$ and the conjugate transpose
over $\C$ or $\HH$. We also write $A^T$ for the real transpose and
denote the space of real skew-symmetric matrices by
$\Skew(n)=\{X\in\R^{n\times n}:X^T=-X\}$.
The notation $Sp(n)$ always refers to the compact quaternionic
unitary group.

\begin{table}[ht]
\caption{Classical matrix groups and their real dimensions.}
\label{tab:lie-groups}
\centering
\small
\begin{tabular}{@{}lll@{}}
\toprule
Group & Lie algebra & Dimension\\
\midrule
$GL(n,\R)$ & $\R^{n\times n}$ & $n^2$\\
$SL(n,\R)$ & $\{X\in\R^{n\times n}:\tr X=0\}$ & $n^2-1$\\
$O(n),\ SO(n)$ & $\Skew(n)$ & $n(n-1)/2$\\
$U(n)$ & $\{X\in\C^{n\times n}:X^*=-X\}$ & $n^2$\\
$SU(n)$ & $\{X\in\mathfrak u(n):\tr X=0\}$ & $n^2-1$\\
$Sp(n)$ & $\{X\in\HH^{n\times n}:X^*=-X\}$ & $n(2n+1)$\\
\bottomrule
\end{tabular}
\end{table}

The Lie algebra is the tangent space at the identity.
Differentiating $Q^*Q=I$ gives $X^*=-X$, and differentiating
$\det Q=1$ gives $\tr X=0$. The dimensions in
Table~\ref{tab:lie-groups} follow by counting entries. In particular,
a quaternionic skew-adjoint matrix has three real parameters in
each diagonal entry and four in each entry above the diagonal.

\subsection{The common Cayley formula}

Let $\mathcal A\subset\mathbb F^{n\times n}$ be a real unital
algebra closed under taking adjoints, and set
$$
G=\{Q\in\mathcal A:Q^*Q=I\},\qquad
\fkg=\{X\in\mathcal A:X^*=-X\}.
$$
For $\mathcal A=\mathbb F^{n\times n}$, these are the orthogonal,
unitary and compact symplectic groups and their Lie algebras,
respectively. Denote by $G^0$ the identity component of $G$.

\begin{definition}
\label{def:cayley-type}
Set $G_+=\{Q\in G:I+Q\text{ is invertible}\}$. The Cayley
transform is
\begin{equation}
\label{eq:cayley-general}
\mathbf C_n:\fkg\longrightarrow G_+,\qquad
\mathbf C_n(X)=(I+X)^{-1}(I-X).
\end{equation}
\end{definition}

We occasionally write $\mathbf C=\mathbf C_n$ when the matrix
size is fixed.

\begin{theorem}
\label{thm:cayley-properties}
The Cayley transform has the following properties.
\begin{enumerate}
\item[(i)] It is a real-analytic diffeomorphism, with inverse
$$
\mathbf C_n^{-1}(Q)=(I+Q)^{-1}(I-Q).
$$
\item[(ii)] For $X,Y\in\fkg$, we have
\begin{equation}
\label{eq:cayley-differential}
d(\mathbf C_n)_X(Y)=-2(I+X)^{-1}Y(I+X)^{-1}.
\end{equation}
In particular, $d(\mathbf C_n)_0=-2\id$, and $\mathbf C_n$
has no critical points.
\item[(iii)] We have $G_+\subset G^0$, and
$\Sigma^0=G^0\setminus G_+$ is closed and has measure zero.
Thus $\mathbf C_n$ is a strict approximate parametrization of $G^0$.
\item[(iv)] The group $G$ has finitely many components
$R_1G^0,\ldots,R_rG^0$. If $\dim_\R\fkg\geq1$, the map
$$
(j,X)\longmapsto R_j\mathbf C_n(X)
$$
is a strict approximate parametrization of $G$ from
$\{1,\ldots,r\}\times\fkg$, with omitted set
$\bigcup_{j=1}^rR_j\Sigma^0$.
\item[(v)] For every $R\in G$ and $X\in\fkg$,
$$
\mathbf C_n(RXR^{-1})=R\mathbf C_n(X)R^{-1}.
$$
\end{enumerate}
\end{theorem}

\begin{proof}
If $B\in\mathcal A$ is invertible, left multiplication by $B$
is an injective real-linear map on $\mathcal A$. Since
$\mathcal A$ is finite-dimensional, it is surjective. Solving
$BC=I$ in $\mathcal A$ gives $B^{-1}\in\mathcal A$.
For $X^*=-X$ and $v\in\mathbb F^n$, we have
$$
|(I+X)v|^2=|v|^2+|Xv|^2,
$$
because $\text{Re}(v^*Xv)=0$. Hence $I+X$ is invertible.
The same argument applies to $I-X$.

Let $Q=(I+X)^{-1}(I-X)$. Since $I+X$, $I-X$ and their
inverses commute, we obtain
$$
\begin{array}{l}
Q^*=(I+X)(I-X)^{-1},\\[3pt]
Q^*Q=I,\\[3pt]
I+Q=2(I+X)^{-1}.
\end{array}
$$
Thus $Q\in G_+$. Conversely, let $Q\in G_+$ and set
$X=(I+Q)^{-1}(I-Q)$. Using $Q^*=Q^{-1}$, we have
$$
X^*=(I-Q^{-1})(I+Q^{-1})^{-1}
=(Q-I)(Q+I)^{-1}=-X.
$$
Moreover,
$$
I+X=2(I+Q)^{-1},\qquad I-X=2Q(I+Q)^{-1}.
$$
These identities prove that the two maps are inverses.
Matrix inversion is real-analytic where defined, also for
quaternionic matrices viewed as real-linear maps. This proves (i).
The chart near $I$ and its left translates also give the
real-analytic structure on $G$.

Writing $\mathbf C_n(X)=2(I+X)^{-1}-I$ and differentiating
gives \eqref{eq:cayley-differential}. Both factors multiplying
$Y$ are invertible, so the differential has zero kernel.
This proves (ii). Substituting
$I\pm RXR^{-1}=R(I\pm X)R^{-1}$ into
\eqref{eq:cayley-general} also proves (v).

For (iii), the path $t\mapsto\mathbf C_n(tX)$ joins $I$ to
$\mathbf C_n(X)$ in $G$. Hence $G_+\subset G^0$.
The real determinant of $I+Q$, viewed as an endomorphism of the
underlying real vector space of $\mathbb F^n$, is a real-analytic
function of $Q$. It is nonzero at $I$, and its zero set in $G^0$
is $\Sigma^0$.

We recall why this zero set has measure zero. A nonzero
real-analytic function on a connected analytic manifold cannot
vanish to infinite order at any point. At each zero, choose a
nonzero partial derivative of least order. Reducing its order
by one gives a derivative which vanishes at that point and has
nonzero differential there. Its regular zero set is a smooth
hypersurface. Countably many coordinate charts and multi-indices
therefore cover the original zero set by smooth hypersurfaces.
Proposition~\ref{prop:submanifold-measure} gives measure zero,
and continuity gives closedness.

Finally, $G$ is a closed subgroup of the compact group of
$\mathbb F$-linear isometries of $\mathbb F^n$. Its components
are open, so there are finitely many of them. Left translation
carries $G_+$ and $\Sigma^0$ to each component. If
$d=\dim_\R\fkg\geq1$, map each copy of $\fkg\cong\R^d$
onto a unit ball by $x\mapsto x/\sqrt{1+|x|^2}$ and translate
the finitely many balls so that they are disjoint. This identifies
the parameter space in (iv) with an open subset of $\R^d$,
as required by Definition~\ref{def:approx-param}.
\end{proof}

\section{Cayley coordinates on Lie groups}
\label{sec:lie-groups}

In this section, we prove Theorems~\ref{thm:main-classical}
and~\ref{thm:main-special-unitary}. We use the Cayley identities
from Section~\ref{sec:cayley} and compute the dimensions of the
omitted sets by prescribing the $-1$ eigenspaces. We then combine
the orthogonal parametrization with the QR decomposition to obtain
coordinates on $GL(n,\R)$ and $SL(n,\R)$.

\subsection{The orthogonal groups}
\label{subsec:On}
\label{subsec:SO}

Set
$$
SO_+(n)=\{Q\in SO(n):I+Q\text{ is invertible}\},\qquad
\Sigma_n^{SO}=SO(n)\setminus SO_+(n).
$$
An orthogonal matrix without eigenvalue $-1$ has determinant one:
its nonreal eigenvalues occur in conjugate pairs with product one,
and its remaining eigenvalues are all equal to one. Hence
$SO_+(n)=\{Q\in O(n):I+Q\text{ is invertible}\}$.

\begin{theorem}
\label{thm:On-cayley}
\label{thm:SOn-cayley}
Let $n\geq2$. The Cayley transform
$$
\mathbf C_n:\Skew(n)\longrightarrow SO_+(n),\qquad
\mathbf C_n(X)=(I+X)^{-1}(I-X),
$$
is a real-analytic diffeomorphism, with inverse
$\mathbf C_n^{-1}(Q)=(I+Q)^{-1}(I-Q)$.
For $1\leq k\leq n$, set
$$
\Sigma_{n,k}^{SO}
=\{Q\in SO(n):\dim_\R\ker(I+Q)=k\}.
$$
This stratum is nonempty precisely when $k$ is even. Every nonempty
stratum is a smooth embedded submanifold, and
$$
\codim_{SO(n)}\Sigma_{n,k}^{SO}=\frac{k(k-1)}2.
$$
Consequently, $\Sigma_n^{SO}$ is closed, has measure zero and has
codimension one. Thus $\mathbf C_n$ is a strict approximate
parametrization of $SO(n)$. For $n=1$, the same formula gives the
global parametrization of $SO(1)=\{1\}$ by $\R^0$.
\end{theorem}

\begin{proof}
The diffeomorphism and its inverse follow from
Theorem~\ref{thm:cayley-properties}. We also have
$$
\det\mathbf C_n(X)=\frac{\det(I-X)}{\det(I+X)}=1,
$$
since $(I-X)^*=I+X$ over $\R$.

We next compute the dimension of each stratum. Let
$Q\in\Sigma_{n,k}^{SO}$ and $E=\ker(I+Q)$. Then $E^\perp$ is
invariant under $Q$, and $B=Q|_{E^\perp}$ has no eigenvalue $-1$.
Thus $\det B=1$ and $1=\det Q=(-1)^k$, so $k$ is even.
Conversely, for even $k$, choose
$E\in\Gr_\R(n,k)$ and $B\in SO_+(E^\perp)$. The matrix
$Q=(-I_E)\oplus B$ then has $-1$ eigenspace exactly $E$.

These choices identify $\Sigma_{n,k}^{SO}$ with a smooth bundle
over $\Gr_\R(n,k)$ with fiber $SO_+(n-k)$. To see that the
identification is smooth in both directions, choose local
orthonormal frames of $E^\perp$. Near a fixed matrix in the stratum,
a small circle enclosing $-1$ and no other eigenvalue defines a
smooth spectral projection onto $E$. This recovers $E$ smoothly,
and the restriction of $Q$ to $E^\perp$ recovers $B$. The same
local description shows that the stratum is embedded. Hence
$$
\dim\Sigma_{n,k}^{SO}
=k(n-k)+\frac{(n-k)(n-k-1)}2
=\frac{n(n-1)}2-\frac{k(k-1)}2.
$$
When $k=n$, the orthogonal complement contributes no parameters.
Each even $k\leq n$ is realized by
$\mathrm{diag}(-I_k,I_{n-k})$. Therefore $\Sigma_n^{SO}$ is a
finite union of smooth strata of positive codimension. It has
measure zero by Proposition~\ref{prop:submanifold-measure} and is
closed because it is defined by $\det(I+Q)=0$. The largest stratum
has $k=2$, giving codimension one.
\end{proof}

To parametrize both components of $O(n)$, we translate the Cayley
chart by a reflection.

\begin{proposition}
\label{prop:On-components}
Let $n\geq2$ and $R=\mathrm{diag}(-1,1,\ldots,1)$. The map
$$
\mathbf G_O:\{0,1\}\times\Skew(n)\longrightarrow O(n),\qquad
\mathbf G_O(\varepsilon,X)=R^\varepsilon\mathbf C_n(X),
$$
is a strict approximate parametrization with image
$SO_+(n)\sqcup RSO_+(n)$. Its omitted set
$$
\widetilde\Sigma_n^O=\Sigma_n^{SO}\sqcup R\Sigma_n^{SO}
$$
has codimension one.
\end{proposition}

\begin{proof}
The two maps are diffeomorphisms onto open subsets of the two
components of $O(n)$. Left translation preserves measure zero and
stratum dimensions, so the omitted set is closed, has measure zero
and has codimension one. Since $n(n-1)/2\geq1$, the two copies of
$\Skew(n)$ are diffeomorphic to two disjoint open balls in
$\R^{n(n-1)/2}$.
\end{proof}

The restriction $n\geq2$ is needed here. The two points of
$O(1)=\{1,-1\}$ cannot be parametrized by an open subset of
$\R^0$, and no nonempty subset of $O(1)$ has measure zero.

\subsection{The unitary group \texorpdfstring{$U(n)$}{U(n)}}
\label{subsec:Un}

Set
$$
\mathfrak u(n)=\{X\in\C^{n\times n}:X^*=-X\},\qquad
U_+(n)=\{Q\in U(n):I+Q\text{ is invertible}\}.
$$
For unitary matrices, the $-1$ eigenspace can have any complex
dimension between one and $n$.

\begin{theorem}
\label{thm:Un-cayley}
Let $n\geq1$. The Cayley transform
$\mathbf C_n:\mathfrak u(n)\to U_+(n)$ is a real-analytic
diffeomorphism, with inverse
$\mathbf C_n^{-1}(Q)=(I+Q)^{-1}(I-Q)$.
For $1\leq k\leq n$, the stratum
$$
\Sigma_{n,k}^U=\{Q\in U(n):\dim_\C\ker(I+Q)=k\}
$$
is a nonempty smooth embedded submanifold of real codimension
$k^2$. Thus $\Sigma_n^U=U(n)\setminus U_+(n)$ is closed, has
measure zero and has codimension one. In particular,
$\mathbf C_n$ is a strict approximate parametrization of $U(n)$.
\end{theorem}

\begin{proof}
The first assertion follows from
Theorem~\ref{thm:cayley-properties}. As in the orthogonal case,
choosing $E=\ker(I+Q)\in\Gr_\C(n,k)$ and
$Q|_{E^\perp}\in U_+(n-k)$ describes the stratum as a smooth
bundle. Its real dimension is
$$
\dim\Sigma_{n,k}^U=2k(n-k)+(n-k)^2=n^2-k^2.
$$
Each stratum contains $\mathrm{diag}(-I_k,I_{n-k})$. Their union
is the closed set defined by $\det(I+Q)=0$ and has measure zero.
The largest stratum has $k=1$, giving codimension one.
\end{proof}

\subsection{The compact symplectic group
\texorpdfstring{$Sp(n)$}{Sp(n)}}
\label{subsec:Spn}

We use the quaternionic realization
$$
Sp(n)=\{Q\in\HH^{n\times n}:Q^*Q=I\},\qquad
\mathfrak{sp}(n)=\{X\in\HH^{n\times n}:X^*=-X\}.
$$
A skew-adjoint quaternionic matrix has three real parameters for
each diagonal entry and four for each entry above the diagonal.
Hence $\dim_\R\mathfrak{sp}(n)=n(2n+1)$.

\begin{theorem}
\label{thm:Spn-cayley}
Let $n\geq1$ and
$Sp_+(n)=\{Q\in Sp(n):I+Q\text{ is invertible}\}$.
The Cayley transform $\mathbf C_n:\mathfrak{sp}(n)\to Sp_+(n)$
is a real-analytic diffeomorphism, with inverse
$\mathbf C_n^{-1}(Q)=(I+Q)^{-1}(I-Q)$.
For $1\leq k\leq n$, the stratum
$$
\Sigma_{n,k}^{Sp}
=\{Q\in Sp(n):\dim_\HH\ker(I+Q)=k\}
$$
is a nonempty smooth embedded submanifold of real codimension
$k(2k+1)$. Thus $\Sigma_n^{Sp}=Sp(n)\setminus Sp_+(n)$ is
closed, has measure zero and has codimension three. In particular,
$\mathbf C_n$ is a strict approximate parametrization of $Sp(n)$.
\end{theorem}

\begin{proof}
The Cayley identities hold over $\HH$ by
Theorem~\ref{thm:cayley-properties}. Let
$E=\ker(I+Q)$ have quaternionic dimension $k$. Then
$Q|_{E^\perp}\in Sp_+(n-k)$, so the stratum is a smooth bundle
over $\Gr_\HH(n,k)$ with fiber $Sp_+(n-k)$. The smoothness
argument is the same as above: in the complex representation,
the spectral projection onto the $-1$ eigenspace has rank $2k$
and preserves the quaternionic structure. Therefore
$$
\dim\Sigma_{n,k}^{Sp}
=4k(n-k)+(n-k)(2(n-k)+1)
=n(2n+1)-k(2k+1).
$$
Each stratum contains $\mathrm{diag}(-I_k,I_{n-k})$. Their union
is closed because singularity of $I+Q$ is a closed condition,
and it has measure zero by the positive codimensions. The largest
stratum has $k=1$, giving codimension three.
\end{proof}

\subsection{The special unitary group
\texorpdfstring{$SU(n)$}{SU(n)}}
\label{subsec:SUn}

Set
$$
\mathfrak{su}(n)=\{Y\in\mathfrak u(n):\tr Y=0\},\qquad
SU_+(n)=SU(n)\cap U_+(n).
$$
If $X\in\mathfrak u(n)$ has eigenvalues
$i\lambda_1,\ldots,i\lambda_n$, then
$$
\det\mathbf C_n(X)
=\exp\left(-2i\sum_{j=1}^n\arctan\lambda_j\right).
$$
The branch of $\arctan$ is the one fixed in the Introduction.
For $n\geq3$, the trace condition does not imply determinant one.
For example, $X=\mathrm{diag}(i,i,-2i,0,\ldots,0)$ satisfies
$\tr X=0$, whereas $\det\mathbf C_n(X)=(3-4i)/5$.
We therefore apply $\mathbf C_n$ to $Y+isI$ and choose $s$ to
satisfy the determinant condition.

\begin{theorem}
\label{thm:SUn-cayley}
Let $n\geq2$ and $K_n=\{k\in\Z:|k|<n/2\}$. For every
$Y\in\mathfrak{su}(n)$ and $k\in K_n$, there is a unique
$s_k(Y)\in\R$ such that
$$
\sum_{j=1}^n\arctan(\mu_j+s_k(Y))=k\pi,
$$
where $\mu_1,\ldots,\mu_n$ are the eigenvalues of $-iY$, counted
with multiplicity. Each $s_k$ is real-analytic, and
$$
\mathbf C_n^{SU}:K_n\times\mathfrak{su}(n)\longrightarrow SU_+(n),
\qquad
\mathbf C_n^{SU}(k,Y)=\mathbf C_n(Y+i s_k(Y)I),
$$
is a real-analytic diffeomorphism. For each $k\in K_n$, the image
of $\{k\}\times\mathfrak{su}(n)$ is one connected component of
$SU_+(n)$, and these images give all the components without
repetition. The map is a strict approximate parametrization of
$SU(n)$, with omitted set
$\Sigma_n^{SU}=SU(n)\setminus SU_+(n)$ of codimension three for
$n=2$ and one for $n\geq3$. For $n=2$, we have $K_2=\{0\}$ and
$s_0=0$.
\end{theorem}

\begin{proof}
We first solve the scalar equation. For fixed $Y$, we have
$$
\begin{array}{l}
\displaystyle\frac{d}{ds}\tr\arctan(-iY+sI)
=\sum_{j=1}^n\frac{1}{1+(\mu_j+s)^2}>0,\\[8pt]
\displaystyle\lim_{s\to\pm\infty}\tr\arctan(-iY+sI)
=\pm\frac{n\pi}{2}.
\end{array}
$$
Hence the sum takes the value $k\pi$ at a unique point precisely
when $k\in K_n$. The function
$(Y,s)\mapsto\tr\arctan(-iY+sI)$ is real-analytic: near any
fixed Hermitian matrix, the holomorphic functional calculus
applies to $\arctan$ on a neighborhood of its real spectrum.
The real-analytic implicit function theorem gives the local
real-analytic dependence of $s_k$ on $Y$, and uniqueness identifies
these local solutions on overlaps. The determinant identity then
gives $\mathbf C_n^{SU}(k,Y)\in SU_+(n)$.

We next describe the inverse. For $Q\in SU_+(n)$, set
$$
\begin{array}{l}
X=(I+Q)^{-1}(I-Q),\\[4pt]
s=\displaystyle\frac{\tr X}{in},\qquad Y=X-isI,\\[6pt]
k=\displaystyle\frac1\pi\tr\arctan(-iX).
\end{array}
$$
By Theorem~\ref{thm:Un-cayley}, $X\in\mathfrak u(n)$, so
$s\in\R$ and $Y\in\mathfrak{su}(n)$. Since $\det Q=1$, the
determinant identity shows that $k\in\Z$. The range of $\arctan$
gives $|k|<n/2$, and hence $k\in K_n$. Uniqueness in the scalar
equation gives $s=s_k(Y)$. Thus the displayed formulas recover
the unique pair $(k,Y)$. The integer $k$ is locally constant in
$Q$, and the remaining inverse formulas are real-analytic. Each
image for fixed $k$ is therefore open and diffeomorphic to the
connected space $\mathfrak{su}(n)$. These images are precisely
the connected components of $SU_+(n)$.

It remains to compute the dimension of $\Sigma_n^{SU}$. Suppose
$E=\ker(I+Q)$ has complex dimension $1\leq r<n$, and put
$m=n-r$. Then the restriction $B=Q|_{E^\perp}$ satisfies
$$
B\in U_+(m),\qquad \det B=(-1)^r.
$$
The determinant map $U(m)\to S^1$ is a submersion, since its
value along $B\exp(itI/m)$ is $e^{it}\det B$. Thus every
nonempty fiber of its restriction to $U_+(m)$ has dimension
$m^2-1$. Choosing $E\in\Gr_\C(n,r)$ and then $B$ describes
the stratum as a smooth bundle, as in the unitary case. Hence
$$
\dim\{Q\in\Sigma_n^{SU}:\dim_\C\ker(I+Q)=r\}
=2r(n-r)+(n-r)^2-1=n^2-r^2-1.
$$
Every nonempty stratum with $r<n$ therefore has codimension $r^2$.

We now check which strata are nonempty. If $m=1$, the fiber is
nonempty precisely when $r$ is even, in which case $B=1$.
If $m\geq2$, we may take $B=I_m$ for even $r$ and
$B=e^{i\pi/m}I_m$ for odd $r$. Neither choice has eigenvalue
$-1$. Finally, $r=n$ gives the single matrix $Q=-I$, which
belongs to $SU(n)$ precisely when $n$ is even and has codimension
$n^2-1$.

For $n=2$, the stratum $r=1$ is empty. Thus
$\Sigma_2^{SU}=\{-I\}$ has codimension three. The eigenvalues
of $-iY$ are $\mu$ and $-\mu$, so $s_0(Y)=0$.
For $n\geq3$, the stratum $r=1$ is nonempty, as shown by
$$
Q=\mathrm{diag}\bigl(-1,e^{i\pi/(n-1)},\ldots,
e^{i\pi/(n-1)}\bigr).
$$
It has the largest dimension, giving codimension one.
In both cases, $\Sigma_n^{SU}$ is closed and is a finite union
of smooth strata of positive codimension, so it has measure zero.
Since $n^2-1>0$, the finitely many copies of $\R^{n^2-1}$ in
the parameter domain are diffeomorphic to disjoint open balls in
the same Euclidean space. This proves the strict parametrization
assertion.
\end{proof}

There are $n$ components of $SU_+(n)$ for odd $n$, and $n-1$
for even $n$. At $Y=0$, we have $s_k(0)=\tan(k\pi/n)$ and
$\mathbf C_n^{SU}(k,0)=e^{-2\pi i k/n}I$.
For $n=1$, the group $SU(1)=\{1\}$ has a global parametrization
by $\R^0$.

\subsection{The general and special linear groups}
\label{subsec:GL}

We use the QR decomposition with positive diagonal entries.
Let $N(b)$ be the unit upper triangular matrix with entries
$b=(b_{ij})_{i<j}\in\R^{n(n-1)/2}$ above the diagonal, and set
$D(a)=\mathrm{diag}(e^{a_1},\ldots,e^{a_n})$ for $a\in\R^n$.
For $A\in GL^+(n,\R)$, write
$$
A=Q(A)D(a(A))N(b(A)),\qquad Q(A)\in SO(n),
$$
and set
$$
GL_+(n,\R)=\{A\in GL^+(n,\R):Q(A)\in SO_+(n)\}.
$$

\begin{proposition}
\label{prop:GLn-param}
Let $n\geq2$. The map
$$
\begin{array}{l}
\mathbf G_{GL}:\Skew(n)\times\R^n\times\R^{n(n-1)/2}
\longrightarrow GL_+(n,\R),\\[4pt]
\mathbf G_{GL}(X,a,b)=\mathbf C_n(X)D(a)N(b)
\end{array}
$$
is a diffeomorphism and a strict approximate parametrization of
$GL^+(n,\R)$. Its omitted set is
$$
\Sigma_{GL}=\{A\in GL^+(n,\R):Q(A)\in\Sigma_n^{SO}\}
$$
and has codimension one. Left translation by a reflection gives
a chart on the negative-determinant component. Together, the two
charts give a strict approximate parametrization of $GL(n,\R)$.
For $n=1$, the maps $a\mapsto e^a$ and $a\mapsto-e^a$
parametrize the two components globally.
\end{proposition}

\begin{proof}
We first verify the QR decomposition and the smooth dependence of
its factors. Apply Gram--Schmidt to the columns $v_1,\ldots,v_n$
of $A$. At step $j$, set
$$
w_j=v_j-\sum_{i<j}\langle q_i,v_j\rangle q_i,\qquad
t_{jj}=\|w_j\|,\qquad q_j=w_j/t_{jj}.
$$
Since the columns are linearly independent, $t_{jj}>0$.
Together with $t_{ij}=\langle q_i,v_j\rangle$ for $i<j$,
these formulas give $A=QT$, where $Q\in SO(n)$ and $T$ is
upper triangular with positive diagonal entries. All factors
depend smoothly on $A$.

If $Q_1T_1=Q_2T_2$, then $Q_2^*Q_1=T_2T_1^{-1}$ is
orthogonal and upper triangular with positive diagonal entries,
so it is the identity. Hence the decomposition is unique. The
parameters in $T=D(a)N(b)$ are
$$
a_i=\log t_{ii},\qquad b_{ij}=t_{ij}/t_{ii}\quad(i<j).
$$
It follows that multiplication is a diffeomorphism
$$
SO(n)\times\R^n\times\R^{n(n-1)/2}
\longrightarrow GL^+(n,\R).
$$
Restricting the first factor to $SO_+(n)$ and applying
Theorem~\ref{thm:SOn-cayley} gives $\mathbf G_{GL}$. Its inverse
is completed by $X=(I+Q(A))^{-1}(I-Q(A))$.
In the QR product coordinates, the omitted set is
$$
\Sigma_n^{SO}\times\R^n\times\R^{n(n-1)/2}.
$$
It is closed, has measure zero and has codimension one by
Theorem~\ref{thm:product}. Left translation gives the chart on
the other component. Each parameter space is $\R^{n^2}$, and
their disjoint union is diffeomorphic to two disjoint open balls
in $\R^{n^2}$.
\end{proof}

For $SL(n,\R)$, the determinant condition becomes
$\sum_{i=1}^n a_i=0$.

\begin{proposition}
\label{prop:SLn-param}
Let $n\geq2$ and $SL_+(n,\R)=SL(n,\R)\cap GL_+(n,\R)$.
The map
$$
\mathbf G_{SL}:\Skew(n)\times\R^{n-1}\times\R^{n(n-1)/2}
\longrightarrow SL_+(n,\R),
$$
$$
\mathbf G_{SL}(X,t,b)=\mathbf C_n(X)
\mathrm{diag}\left(e^{t_1},\ldots,e^{t_{n-1}},
e^{-\sum_{j=1}^{n-1}t_j}\right)N(b)
$$
is a diffeomorphism and a strict approximate parametrization of
$SL(n,\R)$ by $\R^{n^2-1}$. Its omitted set is
$\{A\in SL(n,\R):Q(A)\in\Sigma_n^{SO}\}$ and has codimension
one. For $n=1$, $SL(1,\R)=\{1\}$ has a global parametrization
by $\R^0$.
\end{proposition}

\begin{proof}
Since $\det(QD(a)N(b))=e^{\sum_i a_i}$, the QR diffeomorphism
restricts to
$$
SL(n,\R)\cong SO(n)\times\R^{n-1}\times\R^{n(n-1)/2},
$$
where $a_n=-\sum_{j<n}a_j$. Apply $\mathbf C_n$ to the first
factor. In these coordinates, the omitted set is
$\Sigma_n^{SO}\times\R^{n-1}\times\R^{n(n-1)/2}$, so the
remaining assertions follow from Theorem~\ref{thm:product}.
\end{proof}

The parameter spaces and exceptional codimensions are collected
in Table~\ref{tab:lie-param}. All dimensions are real.

\begin{table}[ht]
\caption{Strict approximate parametrizations of classical groups.
Here $d=n(n-1)/2$. The orthogonal and linear rows assume $n\geq2$;
$U(n)$ and $Sp(n)$ allow $n\geq1$.}
\label{tab:lie-param}
\renewcommand{\arraystretch}{1.2}
\centering
\begin{tabular}{@{}llll@{}}
\toprule
Group & Parameter space & Omitted set & Codimension\\
\midrule
$O(n)$ & $\{0,1\}\times\R^d$ & $\widetilde\Sigma_n^O$ & $1$\\
$SO(n)$ & $\R^d$ & $\Sigma_n^{SO}$ & $1$\\
$U(n)$ & $\R^{n^2}$ & $\Sigma_n^U$ & $1$\\
$SU(2)$ & $\R^3$ & $\{-I\}$ & $3$\\
$SU(n)$, $n\geq3$ & $K_n\times\R^{n^2-1}$ & $\Sigma_n^{SU}$ & $1$\\
$Sp(n)$ & $\R^{n(2n+1)}$ & $\Sigma_n^{Sp}$ & $3$\\
$GL^+(n,\R)$ & $\R^{n^2}$ & $Q(A)\in\Sigma_n^{SO}$ & $1$\\
$SL(n,\R)$ & $\R^{n^2-1}$ & $Q(A)\in\Sigma_n^{SO}$ & $1$\\
\bottomrule
\end{tabular}
\end{table}

\section{Coordinates on classical homogeneous spaces}
\label{sec:symmetric-spaces}

We now construct the parametrizations in
Theorem~\ref{thm:main-homogeneous} and describe their inverses
and omitted sets. We first consider orthonormal frames, then
pass to subspaces and nested sequences of subspaces. Although
these manifolds are quotients of matrix groups, the group
parametrizations do not by themselves give coordinates on the
quotients. We therefore verify each inverse directly. The
conditions under which the inverse is defined also describe
the omitted set.

We then apply graph coordinates to projective spaces and
Lagrangian planes. Finally, we consider non-compact symmetric
spaces, for which the omitted set can be empty. As in
Section~\ref{sec:theory}, all dimensions and codimensions are
real. For a finite union of smooth strata, the dimension is the
largest stratum dimension. In particular, a finite union of
strata of positive codimension has measure zero and empty
interior.

\subsection{Orthonormal frames}
\label{subsec:stiefel}

We first construct the coordinates on real orthonormal frames.
For $1\leq p\leq n$, write
$$
\St(n,p)=\St_\R(n,p)
=\{Q\in\R^{n\times p}:Q^TQ=I_p\}.
$$
The differential of $Q\mapsto Q^TQ$ at a frame $Q$ is
$H\mapsto Q^TH+H^TQ$. It is onto the symmetric matrices,
since $H=QS/2$ maps to $S$. Hence
$\dim\St(n,p)=np-p(p+1)/2$. This manifold is compact:
it is closed, and every column has norm one. By completing a
frame to an orthogonal matrix, we obtain
\begin{equation}\label{eq:stiefel-quotient}
\St(n,p)\cong O(n)/O(n-p).
\end{equation}
The stabilizer acts on the last $n-p$ coordinates, and local
orthonormal completions give local trivializations of the quotient.
To obtain a strict approximate parametrization, we now give a
map and recover its parameters from the frame itself.

Assume first that $p<n$. Fix the standard frame
$E_0=\begin{pmatrix}I_p\\0\end{pmatrix}$ and write
$Q=\begin{pmatrix}Q_1\\Q_2\end{pmatrix}$, where
$Q_1\in\R^{p\times p}$. As in the Introduction, the image
will consist of the frames for which $I_p+Q_1$ is invertible.
We use the parameter space
$$
\cS_{\St}(n,p)=\Skew(p)\oplus\R^{(n-p)\times p}.
$$
Its dimension is $p(p-1)/2+(n-p)p=\dim\St(n,p)$.

\begin{theorem}\label{thm:stiefel}
Let $1\leq p<n$ and set
$$
\St_+(n,p)=\{Q\in\St(n,p):\det(I_p+Q_1)\ne0\}.
$$
The map
\begin{equation}\label{eq:stiefel-map}
\begin{split}
\mathbf G_{\St}:\cS_{\St}(n,p)&\longrightarrow\St_+(n,p),\\
(A,B)&\longmapsto
\begin{pmatrix}
2(I_p+A+B^TB)^{-1}-I_p\\
2B(I_p+A+B^TB)^{-1}
\end{pmatrix}
\end{split}
\end{equation}
is a diffeomorphism, with inverse
\begin{equation}\label{eq:stiefel-inverse}
\mathbf G_{\St}^{-1}(Q)=
\left((I_p+Q_1)^{-1}-(I_p+Q_1)^{-T},\,
Q_2(I_p+Q_1)^{-1}\right).
\end{equation}
The omitted set $\Sigma_{\St}=\St(n,p)\setminus\St_+(n,p)$
is closed and real-analytic. Its stratum with
$\dim\ker(I_p+Q_1)=r$ has codimension
$$
r(n-p)+\frac{r(r-1)}2,\qquad 1\leq r\leq p.
$$
Consequently, $\codim\Sigma_{\St}=n-p$, and
$\mathbf G_{\St}$ is a strict approximate parametrization.
\end{theorem}

\begin{proof}
We first check that the map takes values in $\St_+(n,p)$.
Set $M=I_p+A+B^TB$. For $v\ne0$, we have
$v^TMv=|v|^2+|Bv|^2>0$, so $M$ is invertible. Using
$M+M^T=2(I_p+B^TB)$, we obtain
$$
\begin{aligned}
Q_1^TQ_1+Q_2^TQ_2
&=I_p-2(M^{-T}+M^{-1})+4M^{-T}(I_p+B^TB)M^{-1}\\
&=I_p.
\end{aligned}
$$
Since $I_p+Q_1=2M^{-1}$, the image lies in $\St_+(n,p)$.

Conversely, let $Q\in\St_+(n,p)$ and set
$H=(I_p+Q_1)^{-1}$ and $B=Q_2H$. Since $Q_1H=I_p-H$
and $Q_1^TQ_1+Q_2^TQ_2=I_p$, we have
$$
B^TB=H^TH-(I_p-H)^T(I_p-H)=H+H^T-I_p.
$$
Thus $A=H-H^T$ is skew-symmetric and
$I_p+A+B^TB=2H$. Substituting into \eqref{eq:stiefel-map}
gives the original frame $Q$. If instead we start with $(A,B)$,
then $H=M/2$, so $H-H^T=A$ and $Q_2H=B$. This proves
both inverse relations. Matrix inversion is smooth where defined,
and hence the map and its inverse are smooth.

It remains to describe the omitted set and compute its codimension.
Let $E=\ker(I_p+Q_1)$ have dimension $r$. For $v\in E$,
we have $Q_1v=-v$, and hence
$$
|Q_2v|^2=|v|^2-|Q_1v|^2=0.
$$
Thus $Qv=(-v,0)^T$, and orthogonality gives
$Q(E^\perp)\subset(E\oplus\{0\})^\perp$. After choosing $E$,
we therefore choose a $(p-r)$-frame in a space of dimension
$n-r$. Denote its top block, relative to $E^\perp\oplus\R^{n-p}$,
by $\widetilde Q_1$. The kernel is exactly $E$ if and only if
$\det(I_{p-r}+\widetilde Q_1)\ne0$.

On a set of constant rank, the kernel depends smoothly on $Q$.
Using local orthonormal bases, we identify this stratum with a
smooth bundle over the real Grassmannian $\Gr(p,r)$, with fiber
the indicated open subset of $\St(n-r,p-r)$. For $r=p$,
the fiber is a point. The dimension of the stratum is therefore
$$
r(p-r)+(n-r)(p-r)-\frac{(p-r)(p-r+1)}2.
$$
Subtracting from $\dim\St(n,p)=np-p(p+1)/2$ gives the
stated codimension. Every stratum is nonempty: choose $Q=-I$
on $E$ and the standard inclusion on $E^\perp$. The least
codimension is $n-p$, attained when $r=1$.

The equation $\det(I_p+Q_1)=0$ shows that $\Sigma_{\St}$
is closed and real-analytic. Since $p<n$, all its strata have
positive codimension. Proposition~\ref{prop:submanifold-measure}
therefore gives measure zero. Together with the inverse formula,
this proves that $\mathbf G_{\St}$ is a strict approximate
parametrization.
\end{proof}

\begin{remark}\label{rem:stiefel-deriv}
We can also obtain \eqref{eq:stiefel-map} directly from the
Cayley transform. Set
$$
K=\begin{pmatrix}A&B^T\\-B&0\end{pmatrix}\in\Skew(n).
$$
Solving $(I_n+K)Y=E_0$ gives
$Y=(M^{-1},BM^{-1})^T$, where $M=I_p+A+B^TB$. Hence
$$
\mathbf C_n(K)E_0
=(I_n+K)^{-1}(I_n-K)E_0=2Y-E_0=\mathbf G_{\St}(A,B).
$$
Thus $\mathbf G_{\St}(A,B)$ consists of the first $p$ columns
of $\mathbf C_n(K)$. Formula~\eqref{eq:stiefel-inverse} shows
that these columns determine $A$ and $B$ uniquely.
\end{remark}

We next record the square case and two endpoint cases. The square
case requires both components of $O(n)$, while the one-frame case
recovers stereographic coordinates.

\begin{corollary}\label{cor:stiefel-special}
The following hold.
\begin{enumerate}
\item[(i)] If $p=n$, the formula reduces to
$$
A\longmapsto\mathbf C_n(A)=(I_n+A)^{-1}(I_n-A).
$$
Its image is the subset of $SO(n)$ without the eigenvalue $-1$.
It is not dense in $\St(n,n)=O(n)$; multiplication by a fixed
reflection gives a second chart on the other component.
\item[(ii)] If $p=1<n$, the formula is
$$
b\longmapsto
\begin{pmatrix}(1-|b|^2)/(1+|b|^2)\\2b/(1+|b|^2)\end{pmatrix},
\qquad b=\frac{Q_2}{1+Q_1}.
$$
It parametrizes $S^{n-1}\setminus\{(-1,0,\ldots,0)\}$.
The omitted point has codimension $n-1$.
\item[(iii)] If $p=n-1$ and $n\geq2$, positively oriented
completion identifies $\St(n,n-1)$ with $SO(n)$. Under this
identification, Theorem~\ref{thm:stiefel} gives a strict approximate
parametrization whose omitted set has codimension one.
\end{enumerate}
\end{corollary}

\begin{proof}
For (i), set $B=0$ and use the Cayley inverse. A real orthogonal
matrix without eigenvalue $-1$ has determinant one, since its
nonreal eigenvalues occur in conjugate pairs and its real
eigenvalues are $1$. Translation by a reflection gives the chart
on the other component, as in Section~\ref{sec:lie-groups}.
For (ii), substitute $\Skew(1)=\{0\}$ into the map and its
inverse. For (iii), an $(n-1)$-frame has a unique last column
making the completed basis positively oriented. This completion
is smooth, and the codimension is $n-p=1$.
\end{proof}

The same construction applies to complex and quaternionic frames.
Over $\HH$, we use right vector spaces, with matrices acting
on the left. As in Section~\ref{sec:cayley}, $*$ denotes the
conjugate transpose. We retain the order of all matrix factors.
Write
$$
\begin{array}{l}
\St_{\mathbb F}(n,p)=\{Q\in\mathbb F^{n\times p}:Q^*Q=I_p\},\\[3pt]
\mathfrak u(p)=\{A\in\C^{p\times p}:A^*=-A\},\\[3pt]
\mathfrak{sp}(p)=\{A\in\HH^{p\times p}:A^*=-A\}.
\end{array}
$$
Completing a frame to a unitary or quaternionic unitary matrix
gives $U(n)/U(n-p)$ and $Sp(n)/Sp(n-p)$, respectively.
The square cases are included in the following proposition.

\begin{proposition}\label{prop:stiefel-complex}
Let $\mathbb F\in\{\C,\HH\}$ and $1\leq p\leq n$. The map
$$
(A,B)\longmapsto
\begin{pmatrix}2M^{-1}-I_p\\2BM^{-1}\end{pmatrix},
\qquad M=I_p+A+B^*B,
$$
is a diffeomorphism from
$$
\begin{cases}
\mathfrak u(p)\oplus\C^{(n-p)\times p},&\mathbb F=\C,\\
\mathfrak{sp}(p)\oplus\HH^{(n-p)\times p},&\mathbb F=\HH,
\end{cases}
$$
onto the set of frames for which $I_p+Q_1$ is invertible. Its inverse is
$$
A=(I_p+Q_1)^{-1}-(I_p+Q_1)^{-*},\qquad
B=Q_2(I_p+Q_1)^{-1}.
$$
The real dimensions and the codimensions of the omitted strata are
$$
\begin{array}{c@{\quad}|@{\quad}c@{\quad}|@{\quad}c}
\mathbb F&\dim_\R\St_{\mathbb F}(n,p)
&\codim\{\dim_{\mathbb F}\ker(I_p+Q_1)=r\}\\ \hline
\C&2np-p^2&2r(n-p)+r^2\\
\HH&4np-2p^2+p&4r(n-p)+2r^2+r
\end{array}
$$
for $1\leq r\leq p$. The omitted set is closed and null, with
codimension $2(n-p)+1$ over $\C$ and $4(n-p)+3$ over $\HH$.
Both maps are strict approximate parametrizations, including when $p=n$.
\end{proposition}

\begin{proof}
We first check the dimensions. The Hermitian matrices have real
dimension $p^2$ over $\C$ and $2p^2-p$ over $\HH$.
The differential $H\mapsto Q^*H+H^*Q$ is onto, since
$H=QS/2$ maps to $S$. This gives the stated dimensions of
the Stiefel manifolds. The skew-adjoint matrices have real
dimensions $p^2$ and $p(2p+1)$, respectively. Indeed, a
quaternionic skew-adjoint matrix has three real parameters in
each diagonal entry and four in each entry above the diagonal.
Adding the parameters in $B$ gives the same dimensions as the
corresponding manifolds.

For $v\ne0$, we have
$\text{Re}(v^*Mv)=|v|^2+|Bv|^2>0$. Hence $M$ is invertible.
The identity $M+M^*=2(I_p+B^*B)$ gives $Q^*Q=I_p$ by the
calculation in Theorem~\ref{thm:stiefel}, with $*$ replacing $T$.
For the inverse, set $H=(I_p+Q_1)^{-1}$ and $B=Q_2H$.
Preserving the order of the matrix factors, we obtain
$$
B^*B=H^*H-(I_p-H)^*(I_p-H)=H+H^*-I_p.
$$
Thus $A=H-H^*$ gives $I_p+A+B^*B=2H$. The two inverse
relations now follow as in Theorem~\ref{thm:stiefel}.
All these operations are smooth. The calculation remains valid
over $\HH$, since it uses associativity and the adjoint identity
with the factors in the displayed order.

We next compute the omitted codimensions. Let
$E=\ker(I_p+Q_1)$ have $\mathbb F$-dimension $r$.
As in the real case, $Q$ is $-I$ on $E$. Its restriction to
$E^\perp$ is an orthonormal $(p-r)$-frame in
$\mathbb F^{n-r}$. Write $\widetilde Q_1$ for its top block.
The kernel is exactly $E$ when $I_{p-r}+\widetilde Q_1$ is
invertible. Thus each stratum is a smooth
bundle over $\Gr_{\mathbb F}(p,r)$, of real dimension
$$
\begin{aligned}
&2r(p-r)+2(n-r)(p-r)-(p-r)^2, &&\mathbb F=\C,\\
&4r(p-r)+4(n-r)(p-r)-2(p-r)^2+(p-r), &&\mathbb F=\HH.
\end{aligned}
$$
Here the dimension of $\Gr_{\mathbb F}(p,r)$ is
$(\dim_\R\mathbb F)r(p-r)$, as follows by writing nearby
$r$-planes as graphs. Subtracting these dimensions from the
ambient dimensions gives the stated codimensions. Every stratum
occurs by taking $Q=-I$ on $E$ and the standard inclusion
on $E^\perp$. The union is closed, since it is defined by a
rank condition. The least codimensions occur at $r=1$ and
are positive, also when $p=n$. Hence the omitted sets have
measure zero by Proposition~\ref{prop:submanifold-measure}.
This proves the strict approximate parametrization assertions.
\end{proof}

We also record the connectivity used below for the Grassmannians.
The map which forgets the last vector of a frame is a sphere
bundle, with fiber $S^{n-p}$, $S^{2(n-p)+1}$ or
$S^{4(n-p)+3}$ over $\R$, $\C$ or $\HH$, respectively.
The one-frame spaces are spheres. Applying the homotopy sequence
inductively, we obtain
$$
\begin{array}{c@{\quad}|@{\quad}c}
\text{Space}&\text{Connectivity}\\ \hline
\St(n,p),\quad p<n&n-p-1\\
\St_\C(n,p)&2(n-p)\\
\St_\HH(n,p)&4(n-p)+2
\end{array}
$$
In particular, $\St(n,p)$ is simply connected for $n-p\geq2$.
The bundle and homotopy-sequence arguments are described in
\cite[Theorem~4.41 and Examples~4.53--4.54, pp.~376, 381--383]{Hat02}.

\subsection{Linear subspaces and graph coordinates}
\label{subsec:grassmannian}

We next prove the Grassmannian part of
Theorem~\ref{thm:main-homogeneous}. Write
$\Gr(n,p)=\Gr_\R(n,p)$ for the manifold of unoriented real
$p$-planes. Passing from a frame to its column space gives the
principal $O(p)$-bundle $\St(n,p)\to\Gr(n,p)$ and, for
$1\leq p<n$, the quotient description
\begin{equation}\label{eq:grassmann-quotient}
\Gr(n,p)\cong O(n)/(O(p)\times O(n-p)).
\end{equation}
To obtain coordinates and their inverse, write $\R^n=E\oplus F$,
where $E=\R^p\times\{0\}$ and
$F=\{0\}\times\R^{n-p}$. If $V$ is transverse to $F$,
then the projection $V\to E$ is an isomorphism. Thus $V$
is the graph of a unique linear map $B:E\to F$.
The graph itself determines the parameter $B$, independently
of the choice of a frame.

\begin{theorem}\label{thm:grassmannian}
Let $1\leq p<n$ and set
$\Gr_+(n,p)=\{V\in\Gr(n,p):V\cap F=\{0\}\}$. The map
\begin{equation}\label{eq:grassmann-map}
\begin{split}
\mathbf G_{\Gr}:\R^{(n-p)\times p}&\longrightarrow\Gr_+(n,p),\\
B&\longmapsto\operatorname{colspan}
\begin{pmatrix}(I_p+B^TB)^{-1/2}\\B(I_p+B^TB)^{-1/2}\end{pmatrix}
\end{split}
\end{equation}
is a diffeomorphism. For any orthonormal frame
$Q=(Q_1^T,Q_2^T)^T$ of a plane in the image, the inverse is
$B=Q_2Q_1^{-1}$.
The corresponding orthogonal projection is
\begin{equation}\label{eq:grassmann-proj}
P(B)=\begin{pmatrix}I_p\\B\end{pmatrix}
(I_p+B^TB)^{-1}\begin{pmatrix}I_p&B^T\end{pmatrix}.
\end{equation}
The omitted set
$$
\Sigma_{\Gr}=\{V:\dim(V\cap F)\geq1\}
$$
is closed and real-algebraic. Its stratum with $\dim(V\cap F)=r$
has codimension $r^2$, for $1\leq r\leq\min(p,n-p)$.
Thus $\codim\Sigma_{\Gr}=1$, and the map is a strict approximate
parametrization of a manifold of dimension $p(n-p)$.
\end{theorem}

\begin{proof}
We first verify the map and its inverse. Let
$H=(I_p+B^TB)^{-1/2}$, where the square root is positive
definite. Then
$$
\begin{pmatrix}H\\BH\end{pmatrix}^{T}
\begin{pmatrix}H\\BH\end{pmatrix}
=H(I_p+B^TB)H=I_p.
$$
Hence the displayed matrix is an orthonormal frame for the graph
of $B$. Conversely, every $V$ transverse to $F$ is such a
graph, and any frame for $V$ satisfies $Q_2=BQ_1$. This
gives $B=Q_2Q_1^{-1}$. Replacing the frame $Q$ by $QR$
does not change the inverse, since
$(Q_2R)(Q_1R)^{-1}=Q_2Q_1^{-1}$.

Multiplying the normalized frame by its transpose gives
\eqref{eq:grassmann-proj}. This expression is smooth in $B$,
and its inverse is $B=P_{21}P_{11}^{-1}$, which is smooth
where $P_{11}$ is invertible. The graph charts about the
other coordinate planes give the manifold structure and the
dimension $p(n-p)$.

We now describe the planes omitted by the graph map. Let
$L=V\cap F$ and $W=\operatorname{pr}_E(V)$, with
$\dim L=r$ and $\dim W=p-r$. The quotient $V/L$ is the
graph of a unique map $W\to F/L$. Conversely, after choosing
$L$, $W$ and this map, taking the inverse image of the graph
in $E\oplus F$ recovers a unique $p$-plane with intersection
$L$. The stratum is therefore a vector bundle over
$\Gr(n-p,r)\times\Gr(p,p-r)$ with fiber dimension
$(p-r)(n-p-r)$. Its total dimension is
$$
r(n-p-r)+r(p-r)+(p-r)(n-p-r)=p(n-p)-r^2.
$$
Thus its codimension is $r^2$. Every value
$1\leq r\leq\min(p,n-p)$ occurs, so the least codimension
is one.

For an orthonormal frame, $P_{11}=Q_1Q_1^T$ is invertible
if and only if the projection $V\to E$ is an isomorphism.
Hence $\Sigma_{\Gr}$ is defined in projection coordinates
by $\det P_{11}=0$ and is closed and real-algebraic.
The strata all have positive codimension, so
Proposition~\ref{prop:submanifold-measure} gives measure zero.
The map is therefore a strict approximate parametrization.
\end{proof}

\begin{remark}\label{rem:grassmann-chart}
The omitted set is the Schubert incidence locus determined by $F$.
The stratum with $\dim(V\cap F)=1$ is a smooth hypersurface,
while the higher strata have codimensions $4,9,\ldots$ when
they occur. Their union need not be a submanifold. For example,
near $F\in\Gr(4,2)$, write a plane as
$\{(Cy,y):y\in F\}$. The omitted set is then given by
$\det C=0$ in $\R^{2\times2}$, which is singular at $C=0$.
\end{remark}

\begin{remark}\label{rem:grassmann-plucker}
The graph chart is also an affine patch of the Pl\"ucker embedding.
Given a basis $v_1,\ldots,v_p$ of $V$, the line
$[v_1\wedge\cdots\wedge v_p]$ is independent of the basis.
It determines $V$ through
$v\wedge(v_1\wedge\cdots\wedge v_p)=0$.
Where the first Pl\"ucker coordinate is nonzero, the entries
of $B$ are ratios of Pl\"ucker coordinates, with the signs
specified by the determinant expansion. The same argument on
each coordinate patch gives an immersion with a smooth local
inverse. Since the Grassmannian is compact, the Pl\"ucker map
is a projective embedding.
\end{remark}

We next apply the graph construction over $\C$ and $\HH$.
As in the Introduction, set $\beta=\dim_\R\mathbb F$.
The same choice of an intersection, a projected subspace and
a linear map will compute the omitted codimensions.

\begin{proposition}\label{prop:grassmannian-complex}
Let $\mathbb F\in\{\C,\HH\}$ and $1\leq p<n$. The manifolds
$$
\begin{aligned}
&\Gr_\C(n,p)=U(n)/(U(p)\times U(n-p)),\\
&\Gr_\HH(n,p)=Sp(n)/(Sp(p)\times Sp(n-p))
\end{aligned}
$$
have strict approximate parametrizations by
$\mathbb F^{(n-p)\times p}$, given by
$$
B\longmapsto\operatorname{colspan}_{\mathbb F}
\left[\begin{pmatrix}I_p\\B\end{pmatrix}(I_p+B^*B)^{-1/2}\right].
$$
The image consists of planes transverse to
$F=\{0\}\times\mathbb F^{n-p}$, and the inverse is
$B=Q_2Q_1^{-1}$. The real dimension is $\beta p(n-p)$.
For $1\leq r\leq\min(p,n-p)$, the omitted stratum with
$\dim_{\mathbb F}(V\cap F)=r$ has real codimension
$\beta r^2$. In particular, the omitted
set has codimension two over $\C$ and four over $\HH$.
\end{proposition}

\begin{proof}
The quotient descriptions follow by orthonormal completion.
The graph and inverse formulas follow from the proof of
Theorem~\ref{thm:grassmannian}, with $*$ in place of $T$.
Indeed, $I_p+B^*B$ is positive definite and its positive
square root depends smoothly on the real matrix entries.
Right multiplication of a frame by an invertible matrix
preserves its column space, and
$$
(Q_2R)(Q_1R)^{-1}=Q_2Q_1^{-1}
$$
holds also over $\HH$, with the factors in this order.
Thus the inverse does not depend on the frame. In projection
coordinates, it is $B=P_{21}P_{11}^{-1}$ and is smooth.

For the omitted strata, choose $L=V\cap F$,
$W=\operatorname{pr}_E(V)$ and the map $W\to F/L$ as in
the real case, with $E=\mathbb F^p\times\{0\}$.
Each dimension count is multiplied by $\beta$. Hence the
stratum with $\dim_{\mathbb F}L=r$ has real dimension
$\beta(p(n-p)-r^2)$ and codimension $\beta r^2$.
These strata are smooth and nonempty, and their union is
closed by the rank condition. The least codimension is
$\beta>0$, so the omitted set has measure zero. This proves
the remaining assertions.
\end{proof}

Finally, the frame bundles give the following topological
comparisons. Let $1\leq p<n$ and first consider the real
Grassmannian with $n\geq3$. If $n-p\geq2$, then
$\St(n,p)$ is simply connected, and the homotopy sequence gives
$$
\pi_1(\Gr(n,p))\cong\pi_0(O(p))\cong\Z/2.
$$
For $n-p=1$, use the complement map $V\mapsto V^\perp$
and the case $p=1$ above. When $n=2$, the only case is
$\Gr(2,1)=\R P^1\cong S^1$, whose fundamental group is
$\Z$. Over $\C$ and $\HH$, the connectivity of the frame
spaces gives
$$
\begin{aligned}
&\pi_1(\Gr_\C(n,p))=0,\qquad
\pi_2(\Gr_\C(n,p))\cong\pi_1(U(p))\cong\Z,\\
&\pi_j(\Gr_\HH(n,p))=0\quad(1\leq j\leq3),\qquad
\pi_4(\Gr_\HH(n,p))\cong\pi_3(Sp(p))\cong\Z.
\end{aligned}
$$
Here $\pi_1(U(p))\cong\Z$ follows from $U(1)=S^1$ and the
bundles $U(p-1)\to U(p)\to S^{2p-1}$ for $p\geq2$.
Likewise, $Sp(1)=S^3$ and the bundles
$Sp(p-1)\to Sp(p)\to S^{4p-1}$ show that $Sp(p)$ is
2-connected and $\pi_3(Sp(p))\cong\Z$; see
\cite[Examples~4.53--4.54, pp.~382--383]{Hat02}.
Applying the Hurewicz theorem, we obtain the first nonzero
integral homology groups in positive degree:
$$
H_2(\Gr_\C(n,p);\Z)\cong\Z,\qquad
H_4(\Gr_\HH(n,p);\Z)\cong\Z.
$$
Thus $b_2(\Gr_\C(n,p))=b_4(\Gr_\HH(n,p))=1$; see
\cite[Theorem~4.32, p.~366]{Hat02}.

\subsection{Nested subspaces and block elimination}
\label{subsec:flag}

We next construct coordinates for a nested sequence of real
subspaces. Let $0<p_1<\cdots<p_k<n$. Set $p_0=0$,
$p_{k+1}=n$ and $d_j=p_j-p_{j-1}$. The real partial flag
manifold is
$$
\mathcal F(n;p_1,\ldots,p_k)
=\{V_1\subset\cdots\subset V_k\subset\R^n:\dim V_j=p_j\}.
$$
All subspaces are unoriented. Choosing orthonormal bases of the
successive orthogonal differences gives
$$
\mathcal F(n;p_1,\ldots,p_k)
\cong O(n)/(O(d_1)\times\cdots\times O(d_{k+1})),
$$
and hence
$$
\dim\mathcal F
=\frac{n(n-1)}2-\sum_{j=1}^{k+1}\frac{d_j(d_j-1)}2
=\sum_{j=1}^k d_j(n-p_j).
$$
Forgetting the last subspace gives the bundle
\begin{equation}\label{eq:flag-fibration}
\mathcal F(n;p_1,\ldots,p_k)
\longrightarrow\mathcal F(n;p_1,\ldots,p_{k-1})
\end{equation}
with fiber $\Gr(n-p_{k-1},p_k-p_{k-1})$. Indeed, extending a
fixed flag amounts to choosing a $(p_k-p_{k-1})$-plane in
$V_{k-1}^\perp$. Local orthonormal frames of this complement
give local trivializations. Thus one can choose graph
coordinates successively. We will instead write all parameters
in one block lower triangular matrix and obtain its inverse by
block elimination.

Fix $\R^n=E_1\oplus\cdots\oplus E_{k+1}$, where
$\dim E_j=d_j$. Write
$E_{\leq j}=E_1\oplus\cdots\oplus E_j$ and
$E_{>j}=E_{j+1}\oplus\cdots\oplus E_{k+1}$.

\begin{theorem}\label{thm:flag}
Assume $k\geq1$. The partial flag manifold
$\mathcal F(n;p_1,\ldots,p_k)$ admits a strict approximate
parametrization with parameter space
$$
\cS_{\mathcal F}=\prod_{j=1}^k\R^{(n-p_j)\times d_j}.
$$
Regard these parameters as the blocks $L_{ij}$, $i>j$, of a
block lower triangular matrix $L$ with diagonal blocks
$I_{d_j}$. The map sending $L$ to the flag whose $j$th
subspace is spanned by its first $p_j$ columns is a
diffeomorphism onto
$$
\mathcal F_+
=\{(V_1,\ldots,V_k):V_j\cap E_{>j}=\{0\}\text{ for all }j\}.
$$
This set is open and dense. Its complement is closed, has
measure zero and has codimension one.
\end{theorem}

\begin{proof}
We first check that the map takes values in $\mathcal F_+$.
Every leading block of $L$ is invertible, so its first $p_j$
columns span a plane transverse to $E_{>j}$. These column
spaces are nested and therefore define a flag in
$\mathcal F_+$.

To construct the inverse, let $(V_1,\ldots,V_k)\in\mathcal F_+$
and choose an invertible adapted matrix $R$, whose first
$p_j$ columns span $V_j$. The transversality conditions say
precisely that the leading $p_j\times p_j$ blocks of $R$ are
invertible. We claim that there is a unique factorization
$R=LU$, where $L$ has the form stated above and $U$ is
invertible and block upper triangular.

The first diagonal block of $R$ is invertible, so we can
eliminate the blocks below it. Every leading block of the
remaining Schur complement is also invertible: its determinant
multiplied by the determinant of the first pivot block equals
the determinant of the corresponding leading block of $R$.
Induction gives the factorization. To prove uniqueness, suppose
$LU=\widetilde L\widetilde U$. Then
$$
L^{-1}\widetilde L=U\widetilde U^{-1}.
$$
The left side is block lower triangular with identity diagonal,
whereas the right side is block upper triangular. Hence both
sides equal the identity.

We next check that $L$ depends only on the flag. Any other
adapted matrix has the form $RU_0$, where $U_0$ is invertible
and block upper triangular. Since $RU_0=L(UU_0)$, uniqueness
of the factorization leaves $L$ unchanged. The elimination
formulas use matrix multiplication and inversion of invertible
pivot blocks. Applied to local adapted frames, they therefore
give a smooth inverse. The forward map is smooth in
Grassmannian coordinates, so the two maps are diffeomorphisms.

It remains to describe the omitted set. For each $j$, let
$\pi_j:\mathcal F\to\Gr(n,p_j)$ send a flag to $V_j$.
Local orthonormal completions identify the fibers with flags
in $V_j$ and $V_j^\perp$. Thus $\pi_j$ is a surjective smooth
submersion. We have
$$
\mathcal F\setminus\mathcal F_+
=\bigcup_{j=1}^k\pi_j^{-1}
\{V\in\Gr(n,p_j):V\cap E_{>j}\ne\{0\}\}.
$$
By Theorem~\ref{thm:grassmannian}, each set inside braces has
smooth strata of codimensions $r^2$. Their preimages under the
submersion $\pi_j$ have the same codimensions. Since $\pi_j$
is surjective, a nonempty stratum of codimension one occurs
for every $j$. Hence the finite union has measure zero, empty
interior and dimension $\dim\mathcal F-1$. The rank conditions
also show that it is closed. Finally,
$\dim\cS_{\mathcal F}=\sum_{j=1}^kd_j(n-p_j)=\dim\mathcal F$,
as required.
\end{proof}

\begin{remark}\label{rem:flag-bruhat}
The image $\mathcal F_+$ is the usual open cell determined by
the fixed coordinate flag. Block elimination gives its inverse
in the chosen splitting. Equivalently, using the local
trivializations in \eqref{eq:flag-fibration}, one can obtain
these coordinates by successive Grassmannian graph
constructions. Both descriptions parametrize unoriented flags.
\end{remark}

\begin{corollary}\label{cor:complete-flag}
For $n\geq2$, the complete real flag manifold
$\mathcal F(n)=O(n)/(O(1))^n$ has dimension $n(n-1)/2$ and a
strict approximate parametrization by $\R^{n(n-1)/2}$.
The parametrization is given by lower triangular matrices with
diagonal entries one, and its omitted set has codimension one.
For $n=1$, the flag manifold is a point and the omitted set is
empty.
\end{corollary}

\begin{proof}
For $n\geq2$, take $d_j=1$ in Theorem~\ref{thm:flag}. For
$n=1$, there is no proper nonzero subspace to choose.
\end{proof}

\subsection{Projective spaces and Lagrangian planes}
\label{subsec:complex-analogue}

We first record the case of one-dimensional subspaces. For
$n\geq2$, we have
$$
\begin{array}{l}
\C P^{n-1}=U(n)/(U(1)\times U(n-1))=\Gr_\C(n,1),\\
\HH P^{n-1}=Sp(n)/(Sp(1)\times Sp(n-1))=\Gr_\HH(n,1).
\end{array}
$$
The graph maps are the affine charts $z\mapsto[1:z]$ and
$q\mapsto[1:q]$. They omit $\C P^{n-2}$ and $\HH P^{n-2}$,
respectively, of real codimensions two and four. For a right
quaternionic line, the inverse normalization multiplies its
representing column on the right by the inverse of the first
coordinate. If $n=1$, both projective spaces are points and
the omitted sets are empty.

We next impose the Lagrangian condition on a graph. Equip
$\R^{2n}=\R_x^n\oplus\R_y^n$ with the symplectic form
$$
\omega((x,y),(\xi,\eta))=x^T\eta-y^T\xi.
$$
The manifold of Lagrangian planes is
$\Lambda(n)=U(n)/O(n)$. To see this, identify $\R^{2n}$ with
$\C^n$ and choose an orthonormal real basis of a Lagrangian
plane. The symplectic condition makes the imaginary parts of
its Hermitian inner products vanish, so this basis is a
unitary complex basis. Thus $U(n)$ acts transitively, with
stabilizer $O(n)$ at $\R_x^n$. In particular,
$\dim\Lambda(n)=n(n+1)/2$. Denote by $\mathrm{Sym}(n)$ the
space of real symmetric $n\times n$ matrices.

\begin{proposition}\label{prop:lagrangian}
Let $n\geq1$ and $F=\{0\}\times\R^n$. The map
$$
\mathbf G_\Lambda:\mathrm{Sym}(n)\longrightarrow\Lambda(n),
\qquad
S\longmapsto\operatorname{colspan}
\begin{pmatrix}(I_n+S^2)^{-1/2}\\S(I_n+S^2)^{-1/2}\end{pmatrix}
$$
is a diffeomorphism onto the Lagrangian planes transverse to
$F$. If $Q=(Q_1^T,Q_2^T)^T$ is any orthonormal frame of such
a plane, then the inverse is $S=Q_2Q_1^{-1}$.

For $1\leq r\leq n$, the omitted stratum with
$\dim(V\cap F)=r$ has codimension $r(r+1)/2$. The omitted
set is closed, has measure zero and has codimension one.
Hence $\mathbf G_\Lambda$ is a strict approximate
parametrization.
\end{proposition}

\begin{proof}
A plane transverse to $F$ is the graph of a unique real matrix
$S$. On this graph, we have
$$
\omega((x,Sx),(\xi,S\xi))=x^T(S-S^T)\xi.
$$
Thus the plane is Lagrangian if and only if $S=S^T$.
The normalization of its frame and the smooth inverse are
therefore given by the graph calculation in
Theorem~\ref{thm:grassmannian}.

To compute the omitted strata, fix $E=V\cap F$ with
$\dim E=r$. Since $V$ is Lagrangian and contains $E$, we
have $V\subset E^\omega$. The quotient $E^\omega/E$ is
symplectic of dimension $2(n-r)$, and $V/E$ is Lagrangian
and transverse to $F/E$. Conversely, the inverse image of
any such plane is a unique Lagrangian plane $V$ satisfying
$V\cap F=E$.

Identify $E$ with a subspace of $\R_y^n$. Taking its
orthogonal complements in the $x$ and $y$ coordinates
identifies $E^\omega/E$ with $E_x^\perp\oplus E_y^\perp$.
The choice of $V/E$ is then a symmetric graph map on
$E_x^\perp$. Local orthonormal bases give a smooth bundle
over $\Gr(n,r)$ with fiber $\mathrm{Sym}(n-r)$. Hence the
stratum has dimension
$$
r(n-r)+\frac{(n-r)(n-r+1)}2
=\frac{n(n+1)}2-\frac{r(r+1)}2.
$$
The plane $E_x^\perp\oplus E$ realizes each value of $r$.
The union of the strata is closed because it is defined by
a rank condition. Its least codimension is one, so it has
measure zero and the image of $\mathbf G_\Lambda$ is dense.
This proves the remaining assertions.
\end{proof}

For comparison with the non-compact spaces below, we describe
the fundamental group of $\Lambda(n)$. Consider the map
$$
U(n)/O(n)\longrightarrow S^1,\qquad[U]\longmapsto(\det U)^2.
$$
Its fiber over one is $SU(n)/SO(n)$. Indeed, if $\det U=-1$,
right multiplication by a real reflection gives a representative
with determinant one. For a local choice of argument $t$ on
$S^1$, left multiplication by
$\operatorname{diag}(e^{it/2},1,\ldots,1)$ gives a local
trivialization.

For $n\geq2$, the group $SO(n)$ is connected. Also, $SU(n)$
is simply connected, by $SU(2)\cong S^3$ and the bundles
$SU(n-1)\to SU(n)\to S^{2n-1}$. The homotopy sequence
therefore shows that $SU(n)/SO(n)$ is simply connected.
For $n=1$, this fiber is a point. Applying the homotopy
sequence to the displayed bundle gives
$\pi_1(\Lambda(n))\cong\pi_1(S^1)\cong\Z$.

The symmetric matrix $S$ in
Proposition~\ref{prop:lagrangian} parametrizes an open subset
of $U(n)/O(n)$, whose dimension is $n(n+1)/2$. The
non-compact symmetric space $Sp(2n,\R)/U(n)$ has dimension
$n(n+1)$, where $Sp(2n,\R)$ denotes the group of real
$2n\times2n$ symplectic matrices. Its tangent space at the
base point consists of matrices
$$
\begin{pmatrix}A&B\\B&-A\end{pmatrix},\qquad A^T=A,\quad B^T=B.
$$
Indeed, these are precisely the symmetric matrices satisfying
$X^TJ+JX=0$, where
$J=\left(\begin{smallmatrix}0&I\\-I&0\end{smallmatrix}\right)$.
Thus these tangent coordinates require two independent
symmetric matrices.

\subsection{Global coordinates on non-compact symmetric spaces}
\label{subsec:hyperbolic}

We now consider examples for which the omitted set is empty.
We begin with the upper hyperboloid
$$
\mathbb H^n
=\{(x_0,x')\in\R\times\R^n:-x_0^2+|x'|^2=-1,\ x_0>0\}.
$$
Its defining equation determines $x_0$ uniquely from $x'$.

\begin{proposition}\label{prop:hyperbolic}
Let $n\geq1$. The map
$$
\mathbf G_{\mathbb H}:\R^n\longrightarrow\mathbb H^n,
\qquad v\longmapsto(\sqrt{1+|v|^2},v)
$$
is a global diffeomorphism, with inverse
$(x_0,x')\mapsto x'$. Its omitted set is empty.
\end{proposition}

\begin{proof}
The image lies on the upper hyperboloid, since
$-(1+|v|^2)+|v|^2=-1$. Conversely, the defining equation
and $x_0>0$ give $x_0=\sqrt{1+|x'|^2}$. Hence the two
maps are smooth inverses.
\end{proof}

We next give the homogeneous-space description of this model.
Let $SO_0(n,1)$ be the identity component of the group
preserving $-x_0^2+|x'|^2$. Given $v\in\R^n$, set
$a=\sqrt{1+|v|^2}$. The matrix
$$
\begin{pmatrix}a&v^T\\v&I_n+vv^T/(a+1)\end{pmatrix}
$$
sends $(1,0)$ to $(a,v)$ and preserves the quadratic form.
For the latter assertion, write $D=I_n+vv^T/(a+1)$. Then
$Dv=av$ and $D^2=I_n+vv^T$. Together with
$a^2-|v|^2=1$, these identities verify the Lorentz identity
block by block. Replacing $v$ by $tv$ and taking
$a=\sqrt{1+t^2|v|^2}$ gives a path from the identity.
Thus $SO_0(n,1)$ acts transitively on $\mathbb H^n$.
The stabilizer of $(1,0)$ is $SO(n)$, acting on the last
$n$ coordinates, so $\mathbb H^n=SO_0(n,1)/SO(n)$.

\begin{remark}\label{rem:hyperbolic-contrast}
The coordinates above show that $\mathbb H^n$ is contractible.
On its compact dual $S^n$, stereographic projection omits one
point, of codimension $n$. Compactness prevents a global
diffeomorphism $\R^n\to S^n$ when $n\geq1$. As in
Section~\ref{sec:theory}, the empty omitted set has
codimension $+\infty$.
\end{remark}

For the general construction, we use the global Cartan
decomposition. Let $G$ be a connected semisimple Lie group
and let $\theta$ be a Cartan involution of its Lie algebra
$\fkg$. Write $\fkg=\fkk\oplus\fkp$ for the eigenspace
decomposition, and let $K$ be the connected analytic subgroup
with Lie algebra $\fkk$. Then $K$ is closed and the map
$$
K\times\fkp\longrightarrow G,\qquad(k,X)\longmapsto k\exp X
$$
is a diffeomorphism; see~\cite[Theorem~6.31, p.~362]{Kna23}.
If $G$ has finite center, $K$ is maximal compact. For infinite
center, $K$ need not be compact, but this decomposition still
holds.

\begin{theorem}\label{thm:noncompact-symmetric}
Under these assumptions, the map
$$
\mathbf G_{G/K}:\fkp\longrightarrow G/K,
\qquad X\longmapsto(\exp X)K
$$
is a global diffeomorphism. Thus $G/K$ admits a strict global
parametrization by $\R^{\dim\fkp}$, with empty omitted set.
\end{theorem}

\begin{proof}
Since $G/K$ is the quotient by right multiplication, we first
put the $K$-factor on the right. Group inversion and the
identity $(k^{-1}\exp(-X))^{-1}=\exp X\,k$ give a
diffeomorphism
$$
\Phi:\fkp\times K\longrightarrow G,\qquad
(X,k)\longmapsto\exp X\,k.
$$
It follows that every coset has a representative $\exp X$.
If $\exp X\,K=\exp Y\,K$, uniqueness of the factorization
gives $X=Y$. Hence $\mathbf G_{G/K}$ is bijective.

To check the inverse, let $f:G\to\fkp$ be the first component
of $\Phi^{-1}$. For $g=\exp X\,k$ and $k_0\in K$, we have
$gk_0=\exp X\,(kk_0)$, so $f$ is constant on cosets. Since
$K$ is closed, the quotient map admits local smooth sections.
Composing $f$ with these sections gives the smooth inverse
of $\mathbf G_{G/K}$. The forward map is smooth because both
the exponential map and the quotient map are smooth.
\end{proof}

Table~\ref{tab:non-compact-sym} lists the rank-one examples.
In the complex case, we have $S(U(n)\times U(1))\cong U(n)$
by $A\mapsto\operatorname{diag}(A,(\det A)^{-1})$.
The complex and quaternionic examples allow $n\geq1$; see
\cite[Appendix~C.3, pp.~698, 701]{Kna23} for the group pairs
and ranks. For the Cayley hyperbolic plane, $F_4^{-20}$
denotes the simply connected group with real Lie algebra of
type FII. Its Cartan decomposition has
$\fkk=\mathfrak{so}(9)$, and its real rank is one; see
\cite[Appendix~C.4, p.~717]{Kna23}.
Since $G\cong K\times\fkp$ and $G$ is simply connected,
$K$ is simply connected. Thus $K=\mathrm{Spin}(9)$, and
the quotient has dimension $52-36=16$.

\begin{table}[ht]
\caption{Rank-one symmetric spaces of non-compact type.
Each space has a global Euclidean parametrization.}
\label{tab:non-compact-sym}
\renewcommand{\arraystretch}{1.2}
\begin{center}
\begin{tabular}{@{}llll@{}}
\toprule
Space & $G$ & $K$ & Real dimension\\
\midrule
$\mathbb H^n$ ($n\geq2$) & $SO_0(n,1)$ & $SO(n)$ & $n$\\
$\mathbb H^n_\C$ & $SU(n,1)$ & $S(U(n)\times U(1))$ & $2n$\\
$\mathbb H^n_\HH$ & $Sp(n,1)$ & $Sp(n)\times Sp(1)$ & $4n$\\
$\mathbb H^2_{\mathrm{Cay}}$ & $F_4^{-20}$ & $\mathrm{Spin}(9)$ & $16$\\
\bottomrule
\end{tabular}
\end{center}
\end{table}

For comparison, recall the Lie algebra construction of a
compact dual. The bracket relations
$$
[\fkk,\fkk]\subset\fkk,\qquad
[\fkk,\fkp]\subset\fkp,\qquad
[\fkp,\fkp]\subset\fkk
$$
show that $\fkg^*=\fkk\oplus i\fkp$ is a real Lie algebra
in the complexification. For $X,Y\in\fkp$, we have
$[iX,iY]=-[X,Y]$, and the Killing form is negative definite
on $\fkg^*$. Thus $\fkg^*$ is a compact real form; see
\cite[p.~360]{Kna23}. A compact symmetric space requires a
further choice of global groups and subgroups. The sphere
and projective-space charts above give examples. The Lie
algebra construction alone does not determine coordinates
or an omitted set on an arbitrary compact dual.

We collect the parameter dimensions and the codimensions of
the omitted sets in Table~\ref{tab:symmetric-param}. For
frames, the codimension depends on the number of unused
ambient directions. For Grassmannian graph coordinates, the
stratum with one-dimensional intersection has the least
codimension, namely $1$, $2$ or $4$ according to the scalar
field.

\begin{table}[ht]
\caption{Parameter spaces and omitted-set codimensions.
The Grassmannian rows and the real Stiefel row assume $1\leq p<n$;
the other Stiefel rows allow $p=n$. The flag row assumes at least
one proper nonzero subspace.}
\label{tab:symmetric-param}
\renewcommand{\arraystretch}{1.15}
\begin{center}
\small
\begin{tabular}{@{}llll@{}}
\toprule
Space & Parameter space & Real dimension & $\codim\Sigma$\\
\midrule
$\St(n,p)$ & $\Skew(p)\oplus\R^{(n-p)\times p}$
& $np-p(p+1)/2$ & $n-p$\\
$\St_\C(n,p)$ & $\mathfrak u(p)\oplus\C^{(n-p)\times p}$
& $2np-p^2$ & $2(n-p)+1$\\
$\St_\HH(n,p)$ & $\mathfrak{sp}(p)\oplus\HH^{(n-p)\times p}$
& $4np-2p^2+p$ & $4(n-p)+3$\\
$\Gr(n,p)$ & $\R^{(n-p)\times p}$ & $p(n-p)$ & $1$\\
$\Gr_\C(n,p)$ & $\C^{(n-p)\times p}$ & $2p(n-p)$ & $2$\\
$\Gr_\HH(n,p)$ & $\HH^{(n-p)\times p}$ & $4p(n-p)$ & $4$\\
$\mathcal F$ & $\prod_{j=1}^k\R^{(n-p_j)\times d_j}$
& $\sum_{j=1}^kd_j(n-p_j)$ & $1$\\
$\Lambda(n)$ & $\mathrm{Sym}(n)$ & $n(n+1)/2$ & $1$\\
$\mathbb H^n$ & $\R^n$ & $n$ & $+\infty$\\
$G/K$ & $\fkp$ & $\dim\fkp$ & $+\infty$\\
\bottomrule
\end{tabular}
\end{center}
\end{table}

\section{Integration, optimization and further applications}
\label{sec:discussion}

In Section~\ref{sec:theory}, we showed that a strict approximate
parametrization preserves integrals after inserting the Jacobian
and preserves the infimum of every continuous function. We now
apply these facts to the explicit coordinates constructed above.
We first describe gradients, Hessians and retractions in the
parameter domain. We then compute the singular values of the
orthogonal Cayley differential and use them to obtain the Haar
density. Finally, we consider Gaussian models and the rationality
of the orthogonal Cayley chart.

\subsection{Pullback metrics and optimization}
\label{subsec:opt}

Let $\mathbf G:\cS\to\cM\setminus\Sigma$ be a strict
approximate parametrization, where $\cS$ is an open subset of
$\R^d$. Write $U=\cM\setminus\Sigma$. Given a Riemannian
metric $g$ on $\cM$, let $h=\mathbf G^*g$. Then
$$
h_s(v,w)=g_{\mathbf G(s)}(d\mathbf G_s v,d\mathbf G_s w).
$$
For $f\in C^\infty(\cM)$, set $F=f\circ\mathbf G$.
Proposition~\ref{prop:pullback} gives
$$
d\mathbf G_s\bigl(\nabla_hF(s)\bigr)
=\nabla_g f(\mathbf G(s)).
$$
Thus computing the gradient with respect to $h$ is equivalent
to computing the gradient with respect to $g$. If we use the
Euclidean metric on $\cS$, we must also account for the
coefficients of $h$. The existence of Euclidean coordinates
does not imply that the pullback metric is Euclidean.

We next use the same coordinates to define a local retraction.
For $x=\mathbf G(s)$ and
$\xi\in T_x\cM$, define
\begin{equation}
\label{eq:retraction-step}
R_x(\xi)=\mathbf G\bigl(s+(d\mathbf G_s)^{-1}\xi\bigr)
\end{equation}
whenever the argument belongs to $\cS$. Recall that a local
retraction satisfies $R_x(0)=x$ and $d(R_x)_0=\id$.
It is second order if $t\mapsto R_x(t\xi)$ has zero
covariant acceleration at $t=0$.
See~\cite[Section~2.2, Propositions~2--3]{AM09}.

\begin{proposition}
\label{prop:retraction-quality}
The map in~\eqref{eq:retraction-step} is a smooth local retraction
on $U$. The diffeomorphism $\mathbf G:(\cS,h)\to(U,g)$ preserves
the Hessian and its norm. More precisely,
$$
\operatorname{Hess}_hF(v,w)
=\operatorname{Hess}_g f(d\mathbf Gv,d\mathbf Gw).
$$
The retraction in~\eqref{eq:retraction-step} need not be second order.

For $G=SO(n)$ or $U(n)$ with the invariant Frobenius metric,
the recentered Cayley map
$$
\widehat R_Q(Q\Omega)=Q\mathbf C_n(-\Omega/2),
\qquad \Omega^*=-\Omega,
$$
is a second-order retraction, defined on every tangent space.
\end{proposition}

\begin{proof}
Since $\cS$ is open, the pairs $(x,\xi)$ for which
\eqref{eq:retraction-step} is defined form an open subset of
$TU$ containing the zero section. The smoothness of $\mathbf G$
and its inverse gives the smoothness of $R$. By the chain rule,
$R_x(0)=x$ and $d(R_x)_0=\id$. Moreover, $\mathbf G$ is an
isometry from $(\cS,h)$ to $(U,g)$. It therefore preserves
the Levi--Civita connection, which gives the Hessian identity
and the equality of the corresponding tensor norms.

To see that $R$ need not be second order, consider the Cayley
coordinate $a$ on $SO(2)$. The corresponding angle is
$\phi(a)=-2\arctan a$. Along a straight line in the parameter
domain, we have
$$
\frac{d^2}{dt^2}\phi(a+tv)\bigg|_{t=0}
=\frac{4av^2}{(1+a^2)^2}.
$$
The circle metric has constant coefficients in the angle
coordinate. Hence the covariant acceleration at $t=0$ is
nonzero when $a\neq0$ and $v\neq0$.

For $\widehat R$, note that $I-\Omega/2$ is invertible for
every skew-adjoint $\Omega$. Expanding at $t=0$, we obtain
$$
\mathbf C_n(-t\Omega/2)
=I+t\Omega+\frac{t^2}{2}\Omega^2+O(t^3).
$$
The restriction of the Frobenius inner product
$\langle A,B\rangle=\text{Re}\,\tr(A^*B)$ defines a
bi-invariant metric on both groups. Thus $Qe^{t\Omega}$
is the geodesic with initial velocity $Q\Omega$. Its expansion
agrees with that of $\widehat R_Q(tQ\Omega)$ through order two.
This proves that $\widehat R$ is a second-order retraction.
\end{proof}

We now compare the Euclidean Hessian of $F$ with the Riemannian
Hessian of $f$. For constant vectors $v,w$ in the parameter
domain, differentiating twice gives
$$
D^2F[v,w]
=\operatorname{Hess}_g f(d\mathbf Gv,d\mathbf Gw)
+df\bigl((\nabla d\mathbf G)(v,w)\bigr),
$$
where $\nabla d\mathbf G$ uses the flat connection on $\cS$
and the Levi--Civita connection of $g$. If
$$
\|\operatorname{Hess}_g f\|\leq L,\qquad
\|d\mathbf G\|\leq B_1,\qquad
\|\nabla d\mathbf G\|\leq B_2,\qquad
\|\nabla_g f\|\leq M
$$
on the region under consideration, then
$$
\|D^2F\|\leq LB_1^2+MB_2.
$$
Here we use the Euclidean norm on the parameter domain and the
norm induced by $g$ on the target. Smoothness gives such bounds
on each compact subset of $\cS$. A uniform bound on an
unbounded parameter domain requires a separate argument.

For example, the real Stiefel parametrization gives the objective
$$
(A,B)\longmapsto f\bigl(\mathbf G_{\St}(A,B)\bigr).
$$
A local subspace-fitting step in manifold fitting gives a
second, more concrete example; see~\cite{YSLY23,YX25}. Given
observations $z_1,\ldots,z_N\in\R^n$ and $1\leq p<n$, fix a
center $a\in\R^n$, weights $w_i\geq0$ and a smoothing parameter
$\varepsilon>0$, and consider
$$
\mathcal L(V)=\sum_{i=1}^N w_i\sqrt{\varepsilon^2+\|(I-P_V)(z_i-a)\|^2},
\qquad V\in\Gr(n,p),
$$
where $P_V$ is the orthogonal projector onto $V$. This smooth
residual loss reduces the influence of large residuals relative
to a squared loss. Substituting the projector $P(B)$
in~\eqref{eq:grassmann-proj} gives its exact expression
$$
B\longmapsto \mathcal L\bigl(\mathbf G_{\Gr}(B)\bigr)
=\sum_{i=1}^N w_i\sqrt{\varepsilon^2
+\bigl\|(I-P(B))(z_i-a)\bigr\|^2},
\qquad B\in\R^{(n-p)\times p},
$$
in unconstrained matrix coordinates on the chart image of
Theorem~\ref{thm:grassmannian}. The coordinates therefore apply
to a local geometric estimation step within manifold fitting.
Here the center and the weights, and hence the neighborhood they
select, are fixed; we do not study the statistical accuracy of
the resulting estimate.

For continuous $f$, including $f=\mathcal L$,
Theorem~\ref{thm:measure-equiv} gives
$\inf_{\cS}(f\circ\mathbf G)=\inf_{\cM}f$. As explained in
Section~\ref{sec:theory}, this equality does not imply that a
minimizer belongs to the image of the parametrization. Indeed,
if $x_*$ is an omitted point of an embedded compact manifold,
then $f(x)=\|x-x_*\|^2$ has the unique minimizer $x_*$.
Its pullback has infimum zero but no minimizer, and every
minimizing parameter sequence leaves each compact subset of
$\cS$. Wen and Yin~\cite{WY13} use Cayley updates for
orthogonality-constrained optimization. To prove convergence
of an algorithm, one must also specify assumptions on the
objective and the update rule.

We next compute the differential of the orthogonal Cayley map.
This will distinguish the growth of its inverse from the growth
of its condition number near the omitted set. In the following
proposition, both tangent spaces carry the Frobenius norm,
and linear maps carry the induced operator norm.

\begin{proposition}
\label{prop:singular-conditioning}
Let $n\geq2$, $X\in\Skew(n)$, $Q=\mathbf C_n(X)$, and
$A=(I+X)^{-1}$. Denote the singular values of $A$ by
$a_1\geq\cdots\geq a_n>0$. The singular values of
$d\mathbf C_n(X)$ are
$$
2a_i a_j,\qquad 1\leq i<j\leq n.
$$
In particular,
$$
\kappa(d\mathbf C_n(X))
=\frac{a_1a_2}{a_{n-1}a_n},\qquad
\|(d\mathbf C_n(X))^{-1}\|
=\frac1{2a_{n-1}a_n}.
$$
The inverse differential becomes unbounded as $Q$ approaches the
omitted set. The condition number need not become unbounded; it is
identically one when $n=2$.
\end{proposition}

\begin{proof}
Writing $\mathbf C_n(X)=2(I+X)^{-1}-I$ and differentiating,
we obtain
$$
d\mathbf C_n(X)[H]=-2AHA,\qquad
Q^T d\mathbf C_n(X)[H]=-2A^THA.
$$
Left multiplication by $Q^T$ preserves the Frobenius norm.
Choose a singular value decomposition
$A=U\operatorname{diag}(a_i)V^T$. Orthogonal conjugation of
the input by $U$ and the output by $V$ reduces the differential to
$$
H\longmapsto-2\operatorname{diag}(a_i)
H\operatorname{diag}(a_i).
$$
The matrices $(E_{ij}-E_{ji})/\sqrt2$, $i<j$, form an
orthonormal basis of $\Skew(n)$, and the displayed map
multiplies each of them by $-2a_ia_j$. The formulas for the
singular values, the condition number and the inverse norm follow.

To study the omitted set, use $A=(I+Q)/2$ and $a_i\leq1$.
As $Q$ approaches a matrix with eigenvalue $-1$, we have
$a_n\to0$, and hence $a_{n-1}a_n\to0$. Thus the norm of
the inverse differential tends to infinity. When $n=2$,
there is only one singular value, so the condition number
is identically one.
\end{proof}

The quantity $\sigma_{\min}(I+Q)$ and the norm of the inverse
differential therefore describe the loss of control of inverse
coordinates. Note that $I+X$ is always invertible for real
skew-symmetric $X$; all its eigenvalues have modulus at least
one. One can use translated charts to change coordinates before
reaching the omitted set. Its measure zero does not ensure that
deterministic iterates avoid it. If the law of each iterate is
absolutely continuous with respect to volume, then countable
subadditivity shows that all iterates avoid it almost surely.
Their limit may still lie in the omitted set.

\subsection{Haar measure and harmonic analysis}
\label{subsec:harmonic}

We now use the singular values to compute the Jacobian in
Theorem~\ref{thm:measure-equiv}. Let $dQ$ be normalized Haar
measure on $SO(n)$, where $n\geq2$, and let $dX$ be the
Lebesgue measure induced by the Frobenius inner product on
$\Skew(n)$. Set $d=n(n-1)/2$. With the notation of
Proposition~\ref{prop:singular-conditioning}, we have
$$
J_{\mathbf C_n}(X)
=\prod_{i<j}2a_ia_j
=2^d\prod_{i=1}^n a_i^{n-1}
=2^d\det(I+X)^{-(n-1)}.
$$
The last equality follows from $A=(I+X)^{-1}$ and
$\det(I+X)>0$. The positivity follows by pairing the
eigenvalues of $X$. Since the Frobenius volume on $SO(n)$ is
invariant and the omitted set has measure zero, the
change-of-variables formula gives the normalized density
\begin{equation}
\label{eq:cayley-density}
\rho_n(X)=c_n\det(I+X)^{-(n-1)},\qquad
c_n=\frac{2^d}{\Vol_F(SO(n))}.
\end{equation}
Here $\Vol_F$ denotes the Frobenius volume. In particular,
$$
c_n^{-1}=\int_{\Skew(n)}\det(I+X)^{-(n-1)}\,dX<\infty.
$$
Thus compactness of $SO(n)$ gives integrability of this
density on the whole parameter space.

To integrate over $O(n)$, we use one chart on each component,
as in Section~\ref{sec:lie-groups}. Fix a reflection $J\in O(n)$.
Translation by $J$ interchanges the two components, so each
has normalized Haar measure $1/2$. Hence, for every integrable
function $f$ on $O(n)$, we obtain
\begin{equation}
\label{eq:haar-cayley}
\int_{O(n)}f(Q)\,dQ
=\frac12\int_{\Skew(n)}
\bigl(f(\mathbf C_n(X))+f(J\mathbf C_n(X))\bigr)
\rho_n(X)\,dX.
\end{equation}
Here $dQ$ denotes normalized Haar measure on $O(n)$. Removing
the exceptional set does not change this integral. This
assertion concerns Haar measure; it does not apply to surface
measure on a lower-dimensional stratum.

\begin{example}
\label{ex:haar-moments}
Let $Q$ have normalized Haar distribution on $O(n)$. For
every nonnegative integer $k$, formula~\eqref{eq:haar-cayley}
gives
$$
\mathbf E[\tr(Q)^k]
=\frac12\int_{\Skew(n)}
\left(\tr(\mathbf C_n(X))^k+
\tr(J\mathbf C_n(X))^k\right)\rho_n(X)\,dX.
$$
Since $|\tr Q|\leq n$, these integrals are finite. The
integrands and their weights are rational functions of the
Cayley parameters. To check the normalization in dimension
two, write
$$
X=\begin{pmatrix}0&-x\\x&0\end{pmatrix}.
$$
For $SO(2)$, we have $dX=\sqrt2\,dx$, the angle of
$\mathbf C_2(X)$ is $-2\arctan x$, and
$\Vol_F(SO(2))=2\pi\sqrt2$. Substituting into
\eqref{eq:cayley-density}, we obtain
$$
\rho_2(X)\,dX=\frac{dx}{\pi(1+x^2)}.
$$
This verifies both the normalization and the determinant
exponent in~\eqref{eq:cayley-density} when $n=2$.
\end{example}

More generally, let $G$ be a compact Lie group with normalized
Haar measure $dg$, and let
$\mathbf G:\cS\to U=G\setminus\Sigma$ be a strict
approximate parametrization. Denote the pulled-back Haar
measure by $\rho(s)\,ds$. The change-of-variables formula gives
a unitary isomorphism
$$
L^2(G,dg)\longrightarrow L^2(\cS,\rho(s)\,ds),
\qquad f\longmapsto f\circ\mathbf G.
$$
To see surjectivity, compose a function on $\cS$ with
$\mathbf G^{-1}$ on $U$ and extend it arbitrarily across
$\Sigma$. Since $\Sigma$ has measure zero, this extension
defines a unique $L^2$ class. We can therefore apply the
Peter--Weyl theorem in the parameter domain;
see~\cite[Corollary~4.21, p.~247]{Kna23}. As $\pi$ ranges
over the irreducible unitary representations, the functions
$\sqrt{\dim\pi}\,\pi_{ij}\circ\mathbf G$ form an orthonormal
basis with respect to $\rho(s)\,ds$. These functions need
not be spherical harmonics in the parameter variables.

We next consider a transformation preserving Haar measure.
To express all its iterates in one chart, we restrict to
points whose entire orbit stays in $U$. Let
$T:G\to G$ be a Haar-measure-preserving bijection, with $T$ and
$T^{-1}$ measurable. The set
$$
E=\bigcap_{j\in\Z}T^{-j}(U)
$$
is invariant and has full Haar measure, since it is a
countable intersection of sets of full measure. Hence
$\widetilde T=\mathbf G^{-1}T\mathbf G$ is defined on
$\mathbf G^{-1}(E)$ and preserves $\rho(s)\,ds$.
The map $\mathbf G$ identifies invariant measurable sets
modulo sets of measure zero. Thus $\widetilde T$ is ergodic
if and only if $T$ is. The formula for $\widetilde T$ need
not be defined at every point of $\cS$.

\subsection{Gaussian models and the Fisher metric}
\label{subsec:info}

We now consider the effect of coordinates on a statistical
metric. For a regular family of positive densities
$p_\theta$, the Fisher metric is
$$
g_{ij}(\theta)=\mathbf E_\theta
\left[\partial_i\log p_\theta\,\partial_j\log p_\theta\right],
$$
provided the expectations exist and this matrix is positive definite.
The Gaussian family gives a simple example: it admits
global Euclidean parameters, but the coefficients of the
Fisher metric in these parameters need not be constant.

\begin{example}
\label{ex:spd-manifold}
For centered nonsingular Gaussian distributions with covariance
$P\in\mathrm{SPD}(n)$, the Fisher metric is
$$
g_P(U,V)=\frac12\tr(P^{-1}UP^{-1}V).
$$
The map
$$
\mathrm{Sym}(n)\longrightarrow\mathrm{SPD}(n),
\qquad H\longmapsto e^H,
$$
is a global diffeomorphism. For the full Gaussian family, global
coordinates are $(\mu,H)\mapsto(\mu,e^H)$, and the Fisher metric is
$$
g_{(\mu,P)}((u,U),(v,V))
=u^TP^{-1}v+\frac12\tr(P^{-1}UP^{-1}V).
$$
\end{example}

\begin{proof}
We first compute the metric. For $z\sim N(0,P)$,
differentiating the log density in a covariance direction
$U$ gives the score
$$
S_U(z)=\frac12
\left(z^TP^{-1}UP^{-1}z-\tr(P^{-1}U)\right).
$$
Write $z=P^{1/2}y$, where the coordinates of $y$ are
independent standard normal variables. Using their second
and fourth moments, we have, for symmetric matrices $B,C$,
$$
\mathbf E\bigl[(y^TBy-\tr B)(y^TCy-\tr C)\bigr]
=2\tr(BC).
$$
Taking $B=P^{-1/2}UP^{-1/2}$ and
$C=P^{-1/2}VP^{-1/2}$ gives the covariance metric.
For the full Gaussian family, write $z=x-\mu$. The score
in a mean direction $u$ is $u^TP^{-1}z$, which gives
the term $u^TP^{-1}v$. The cross term with $S_U$ vanishes
because the third centered Gaussian moments are zero.

We next verify that the logarithmic coordinates are global.
By the spectral theorem, every $P\in\mathrm{SPD}(n)$ has
a unique symmetric logarithm. The exponential is smooth.
To prove that the logarithm is smooth, including at matrices
with repeated eigenvalues, we use
$$
\log P=\int_0^\infty
\left(\frac{I}{1+t}-(P+tI)^{-1}\right)\,dt.
$$
On each compact subset of $\mathrm{SPD}(n)$, the smallest
eigenvalue has a positive lower bound. The integrand is
uniformly $O(t^{-2})$ as $t\to\infty$, and every derivative
of positive order contains at least two resolvent factors.
Hence the integral and all its derivatives converge locally
uniformly. This proves that the logarithm is smooth and
that the stated maps are diffeomorphisms.
\end{proof}

To write the Fisher metric in these coordinates, we first
differentiate the exponential. For $K\in\mathrm{Sym}(n)$,
we have
$$
d\exp_H[K]=\int_0^1e^{(1-t)H}Ke^{tH}\,dt.
$$
Choose an orthonormal eigenbasis of $H$, with eigenvalues
$h_1,\ldots,h_n$, and express $K,L$ in this basis.
Substituting into the covariance metric, we obtain
$$
(\exp^*g)_H(K,L)
=\frac12\sum_{i,j=1}^n
\left(\frac{\sinh((h_i-h_j)/2)}{(h_i-h_j)/2}\right)^2
K_{ij}L_{ij},
$$
where the quotient has value $1$ at zero. Indeed, the
$(i,j)$ entry of $e^{-H/2}d\exp_H[K]e^{-H/2}$ is
$K_{ij}$ times this quotient. In directions commuting
with $H$, the pullback is $\frac12\tr(KL)$. In general,
its coefficients depend on $H$.

The map $gO(n)\mapsto gg^T$ also identifies
$GL(n,\R)/O(n)$ with $\mathrm{SPD}(n)$. Scalar matrices
give a Euclidean factor. Under the convention excluding
Euclidean factors, this quotient is therefore not a
symmetric space of noncompact type. Note also that the
full Gaussian metric contains the mean term in addition
to the covariance metric.

For comparison, let
$p_\theta=\exp(\theta\cdot T-\psi(\theta))h$ be a regular
minimal exponential family on its open natural parameter domain,
with differentiation under the integral justified. Then
$$
\eta=\nabla\psi(\theta)=\mathbf E_\theta T,\qquad
g=\nabla^2\psi=\operatorname{Cov}_\theta(T).
$$
By minimality,
$a^Tga=\operatorname{Var}_\theta(a\cdot T)>0$ for
$a\neq0$. Since the natural parameter domain is convex,
$\psi$ is strictly convex and $\nabla\psi$ is injective.
The inverse function theorem then gives smooth expectation
coordinates on the image of $\nabla\psi$. The Legendre
transform relates these coordinates to the natural parameters.
A general manifold chart pulls back the Fisher tensor, but
need not be affine for either dual flat connection. Nor
can its exceptional set in general be identified with the
boundary of the statistical family.

\subsection{The rational orthogonal chart}
\label{subsec:birational}

We finally return to the orthogonal Cayley map. In addition
to being smooth, both this map and its inverse are rational
in the matrix entries. We now use this fact to describe
$SO(n)$ as a rational real algebraic variety.

\begin{proposition}
\label{prop:rationality}
The real algebraic group $SO(n)$ is rational. More precisely,
the Cayley map and its inverse,
$$
\mathbf C_n(X)=(I+X)^{-1}(I-X),\qquad
\mathbf C_n^{-1}(Q)=(I+Q)^{-1}(I-Q),
$$
are inverse regular maps between the Zariski open subsets
$$
\{X^T=-X:\det(I+X)\neq0\},\qquad
\{Q\in SO(n):\det(I+Q)\neq0\}.
$$
Translation by a fixed reflection gives the same conclusion for
the other algebraic component of $O(n)$.
\end{proposition}

\begin{proof}
The identities in Theorem~\ref{thm:On-cayley} involve only
transpose, matrix multiplication and inversion. They are
therefore rational identities wherever the denominators
are nonzero. In particular, the displayed maps are inverse
regular maps, and $\det\mathbf C_n(X)=1$. Both open sets
are nonempty: they contain $0$ and $I$, respectively.

To verify density on the real points, use the stratum
calculation in Section~\ref{sec:lie-groups}. It shows that
the image of $\mathbf C_n$ is dense in $SO(n)$ in the
Euclidean topology. Thus every real polynomial vanishing
on this image vanishes on all of $SO(n)$, and the image
is Zariski dense as a real algebraic set. Since it is
isomorphic to a nonempty open subset of the affine space
$\Skew(n)$, this proves rationality. Multiplication by a
fixed reflection gives an algebraic isomorphism onto the
other component of $O(n)$.
\end{proof}

\begin{remark}
\label{rem:birat-cayley}
For $n\geq2$, the Cayley image is a proper open subset of
the identity component. Its omitted set describes the
boundary of these coordinates; it need not be an exceptional
divisor of a resolution. The two algebraic components of
$O(n)$ cannot together be birational to one irreducible
affine space. More generally, the existence of a smooth
strict approximate parametrization alone gives no
conclusion about algebraic rationality.
\end{remark}

\subsection*{Tool use}
The authors used Claude (Anthropic), a large language model, as an
assistant in preparing this manuscript. It was used for LaTeX
editing and typesetting, for checking grammar, the internal
consistency of the text and the cross-references, for checking
selected computations, and for verifying the bibliographic details
of cited works. The mathematical results, their proofs and the
examples are the work of the authors. The authors reviewed and
verified all changes made with this assistance and take full
responsibility for the content of the paper.


\begin{thebibliography}{MVPT18}

\bibitem[AM12]{AM09}
P.-A. Absil and J.~Malick,
Projection-like retractions on matrix manifolds,
\textit{SIAM J. Optim.} \textbf{22} (2012), no.~1, 135--158.

\bibitem[Arn67]{Arn67}
V.~I.~Arnol'd,
Characteristic class entering in quantization conditions,
\textit{Funct. Anal. Appl.} \textbf{1} (1967), no.~1, 1--13.

\bibitem[Cay46]{Cay46}
A.~Cayley,
Sur quelques propri\'et\'es des d\'eterminants gauches,
\textit{J. Reine Angew. Math.} \textbf{32} (1846), 119--123.

\bibitem[FF16]{FF16}
A.~Fomenko and D.~Fuchs,
\textit{Homotopical Topology},
second edition, Graduate Texts in Mathematics, vol.~273,
Springer, Cham, 2016.

\bibitem[Hat02]{Hat02}
A.~Hatcher,
\textit{Algebraic Topology},
Cambridge University Press, Cambridge, 2002.

\bibitem[IT98]{IT98}
J.-i.~Itoh and M.~Tanaka,
The dimension of a cut locus on a smooth Riemannian manifold,
\textit{Tohoku Math. J.} (2) \textbf{50} (1998), no.~4, 571--575.

\bibitem[JHD20]{JHD20}
M.~Jauch, P.~D.~Hoff and D.~B.~Dunson,
Random orthogonal matrices and the Cayley transform,
\textit{Bernoulli} \textbf{26} (2020), no.~2, 1560--1586.

\bibitem[Kna23]{Kna23}
A.~W. Knapp,
\textit{Lie Groups Beyond an Introduction},
Digital second edition, published by the author, 2023.
\url{https://www.math.stonybrook.edu/~aknapp/download/Beyond2-clickable.pdf}.

\bibitem[KY24]{KY24}
K.~Kume and I.~Yamada,
Generalized left-localized Cayley parametrization for optimization with
orthogonality constraints,
\textit{Optimization} \textbf{73} (2024), no.~4, 1113--1159.

\bibitem[Lee18]{Lee18}
J.~M.~Lee,
\textit{Introduction to Riemannian Manifolds},
second edition, Graduate Texts in Mathematics, vol.~176,
Springer, Cham, 2018.

\bibitem[MVPT18]{MVPT18}
E.~Mac\'ias-Virg\'os, M.~J.~Pereira-S\'aez and D.~Tanr\'e,
Cayley transform on Stiefel manifolds,
\textit{J. Geom. Phys.} \textbf{123} (2018), 53--60.

\bibitem[Mil85]{Miller85}
H.~Miller,
Stable splittings of Stiefel manifolds,
\textit{Topology} \textbf{24} (1985), no.~4, 411--419.

\bibitem[Mil65]{Mil65}
J.~W. Milnor,
\textit{Topology from the Differentiable Viewpoint},
University Press of Virginia, Charlottesville, VA, 1965.

\bibitem[SPG25]{SPG25}
K.~Sch\"afers, M.~Peardon and M.~G\"unther,
A modified Cayley transform for $SU(3)$ molecular dynamics simulations,
\textit{Comput. Phys. Commun.} \textbf{315} (2025), 109740.

\bibitem[Sch22]{Sch22}
S.~Schwede,
Global stable splittings of Stiefel manifolds,
\textit{Doc. Math.} \textbf{27} (2022), 789--845.

\bibitem[Toy48]{Toy48}
H.~T\^oyama,
On Haar measure of some groups,
\textit{Proc. Japan Acad.} \textbf{24} (1948), no.~9, 13--16.

\bibitem[WY13]{WY13}
Z.~Wen and W.~Yin,
A feasible method for optimization with orthogonality constraints,
\textit{Math. Program.} \textbf{142} (2013), nos.~1--2, 397--434.

\bibitem[YSLY23]{YSLY23}
Z.~Yao, J.~Su, B.~Li and S.-T.~Yau,
Manifold fitting,
preprint (2023), \texttt{arXiv:2304.07680}.

\bibitem[YX25]{YX25}
Z.~Yao and Y.~Xia,
Manifold fitting under unbounded noise,
\textit{J. Mach. Learn. Res.} \textbf{26} (2025), no.~45, 1--55.

\end{thebibliography}
\end{document}